\documentclass[11pt,a4paper]{article}
\usepackage{amsfonts}
\usepackage{amsmath, amssymb, amsthm}
\usepackage{array,booktabs,color}
\usepackage[T1]{fontenc}
\usepackage{microtype}
\usepackage[american]{babel}
\usepackage{csquotes}
\usepackage{setspace}
\usepackage{blkarray}
\usepackage{afterpage}
\usepackage{diagbox}
\usepackage{dirtree}
\usepackage{enumitem}
\usepackage{eso-pic}
\usepackage{fancybox}
\usepackage{float}
\usepackage{graphicx}
\usepackage{grffile}
\usepackage[hidelinks]{hyperref}
\usepackage[nameinlink]{cleveref}
\newcommand{\textcite}{\citet}

\newcommand{\tablecite}[1]{\citet{#1}}

\usepackage{listings}
\usepackage{longtable}
\usepackage{makecell}
\usepackage{mathtools}
\usepackage{accents}
\usepackage{multirow}
\usepackage{multicol}
\usepackage[short, nocomma]{optidef}
\usepackage{pgfplots}
\usepackage{pgfplotstable}
\usepackage{pifont}
\usepackage{scalerel}
\usepackage{siunitx}
\usepackage{subcaption}
\usepackage{tabularx}
\usepackage{textcomp}
\usepackage{tikz}
\usepackage{varwidth}
\usepackage{yhmath}

\pgfplotsset{compat=1.16}
\usepgfplotslibrary{fillbetween,groupplots}

\makeatletter
\@ifundefined{qty}{\DeclareRobustCommand{\qty}[3][]{\SI[#1]{#2}{#3}}}{}
\@ifundefined{unit}{}{}
\@ifundefined{byte}{\DeclareSIUnit{\byte}{B}}{}
\makeatother

\usetikzlibrary{svg.path, fit, patterns}
\definecolor{orcidlogocol}{HTML}{A6CE39}
\tikzset{
    orcidlogo/.pic={
        \fill[orcidlogocol] svg{M256,128c0,70.7-57.3,128-128,128C57.3,256,0,198.7,0,128C0,57.3,57.3,0,128,0C198.7,0,256,57.3,256,128z};
        \fill[white] svg{M86.3,186.2H70.9V79.1h15.4v48.4V186.2z}
        svg{M108.9,79.1h41.6c39.6,0,57,28.3,57,53.6c0,27.5-21.5,53.6-56.8,53.6h-41.8V79.1z M124.3,172.4h24.5c34.9,0,42.9-26.5,42.9-39.7c0-21.5-13.7-39.7-43.7-39.7h-23.7V172.4z}
        svg{M88.7,56.8c0,5.5-4.5,10.1-10.1,10.1c-5.6,0-10.1-4.6-10.1-10.1c0-5.6,4.5-10.1,10.1-10.1C84.2,46.7,88.7,51.3,88.7,56.8z};
    }
}

\newcommand\orcidicon[1]{\href{https://orcid.org/#1}{\mbox{\scalerel*{
                \begin{tikzpicture}[yscale=-1,transform shape]
                \pic{orcidlogo};
                \end{tikzpicture}
            }{|}}}}

\usetikzlibrary{shapes.geometric, positioning, calc}
\usetikzlibrary{arrows.meta, backgrounds}

\crefname{subtaski}{subtask}{subtasks}

\newcommand{\N}{\ensuremath{\mathbb{N}}}

\newcommand{\R}{\ensuremath{\mathbb{R}}}
\newcommand{\Z}{\ensuremath{\mathbb{Z}}}
\DeclareMathOperator*{\argmax}{arg\,max}

\newcommand{\DemandSym}{g}                   % overall demand at node
\newcommand{\C}{\ensuremath{C}}              % an individual contingency
\newcommand{\DemandParSym}[0]{\ensuremath{G}}  % dynamic demands parameters
\newcommand{\StatDSym}[0]{\ensuremath{D}}    % static demands
\newcommand{\PotSym}[0]{\ensuremath{\theta}} % potentials
\newcommand{\SusSym}[0]{\ensuremath{B}}      % susceptance ?
\newcommand{\FlowSym}[0]{\ensuremath{x}}     % flow variables
\newcommand{\FlowParSym}[0]{\ensuremath{X}}  % flow parameters
\newcommand{\FlowOverSym}[0]{\ensuremath{o}} % overload flow variables
\newcommand{\Arcs}{A}                        % arcs
\newcommand{\LArcs}{\Arcs_L}                 % line arcs
\newcommand{\SArcs}{\Arcs_S}                 % switch arcs
\newcommand{\subSym}{W}                      % substations ?

\newcommand{\TAL}{Z}
\newcommand{\DemandCostsSym}{c}

\newcommand{\ubar}[1]{\underaccent{\bar}{#1}}
\newcommand{\fixedSym}[1]{\hat{{#1}}}
\newcommand{\etaFixedSym}{\ensuremath{\fixedSym{\eta}}}
\newcommand{\lbsym}[1]{\ubar{{#1}}}
\newcommand{\ubsym}[1]{\bar{{#1}}}

\newcommand{\aC}{\ensuremath{\mathcal{C}}} % all contingencies featured in the current instance

\newcommand{\Demand}[1]{\DemandSym_{#1}}

\newcommand{\DemandFixedSym}[0]{\ensuremath{\fixedSym{\DemandSym}}}
\newcommand{\DemandFixed}[1]{\ensuremath{\DemandFixedSym_{{#1}}}}

\newcommand{\StatD}[1]{\ensuremath{\StatDSym_{{#1}}}}

\newcommand{\Pot}[1]{\ensuremath{\PotSym_{#1}}}
\newcommand{\PotCSym}[0]{\ensuremath{\PotSym^{\C}}}
\newcommand{\PotC}[1]{\ensuremath{\PotCSym_{#1}}}

\newcommand{\PotDifMaxSym}[0]{\ensuremath{\bar{\PotSym}}}
\newcommand{\PotDifMax}[2]{\ensuremath{\PotDifMaxSym_{({#1},{#2})}^{\text{on}}}}

\newcommand{\PotDifMaxOff}[2]{\ensuremath{\bar{\PotSym}_{({#1},{#2})}^{\text{off}}}}

\newcommand{\Sus}[2]{\ensuremath{\SusSym_{({#1},{#2})}}}

\newcommand{\FlowMaxSym}[0]{\ensuremath{\bar{\FlowParSym}}}
\newcommand{\FlowMax}[2]{\ensuremath{\FlowMaxSym_{({#1},{#2})}}}

\newcommand{\FlowMaxOverSym}[0]{\ensuremath{\tilde{\FlowParSym}}}
\newcommand{\FlowMaxOver}[2]{\ensuremath{\FlowMaxOverSym_{({#1},{#2})}}}

\newcommand{\FlowOverC}[3]{\ensuremath{\FlowOverSym_{({#1},{#2})}^{{#3}}}}

\newcommand{\FlowOverCOptSym}[0]{\ensuremath{\FlowOverSym^{\C*}}}
\newcommand{\FlowOverCOpt}[2]{\ensuremath{\FlowOverCOptSym_{({#1},{#2})}}}

\newcommand{\FlowCSym}[0]{\ensuremath{\FlowSym^{\C}}}
\newcommand{\FlowC}[2]{\ensuremath{\FlowCSym_{({#1},{#2})}}}
\newcommand{\FlowCOptSym}[0]{\ensuremath{\FlowSym^{\C*}}}
\newcommand{\FlowCOpt}[2]{\ensuremath{\FlowCOptSym_{({#1},{#2})}}}

\newcommand{\OverCost}[0]{\ensuremath{\DemandCostsSym^{\FlowOverSym}}}

\newcommand{\zSym}[0]{\ensuremath{z}}
\newcommand{\z}[2]{\ensuremath{\zSym_{({#1},{#2})}}}
\newcommand{\zDefault}[2]{\ensuremath{\zSym^\text{default}_{({#1},{#2})}}}
\newcommand{\zFixedSym}[0]{\ensuremath{\fixedSym{\zSym}}}
\newcommand{\zFixed}[2]{\ensuremath{\zFixedSym_{({#1},{#2})}}}

\newcommand{\ArcsOutz}{\Arcs^{\zFixedSym,\C}}

\newcommand{\LArcsCSym}[2]{\LArcs^{{#1},{#2}}}
\newcommand{\LArcsC}[1]{\LArcsCSym{\zFixedSym}{#1}}
\newcommand{\LArcsCC}{\LArcsC{\C}}

\newcommand{\SArcsCSym}[1]{\SArcs^{{#1}}}
\newcommand{\SArcsz}{\SArcsCSym{\zFixedSym}}

\newcommand{\Sub}{\subSym}

\newcommand{\Subs}{\mathcal{\subSym}}

\newcommand{\DualFCSym}[0]{\ensuremath{\varphi}}
\newcommand{\DualFCCSym}[0]{\ensuremath{\DualFCSym^\C}}
\newcommand{\DualFCCOptSym}[0]{\ensuremath{\DualFCSym}^{\C*}}
\newcommand{\DualFCCOpt}[1]{\ensuremath{\DualFCCOptSym_{#1}}}
\newcommand{\DualFCC}[1]{\ensuremath{\DualFCCSym_{#1}}}

\newcommand{\DualPESym}[0]{\ensuremath{\omega}}

\newcommand{\DualPECOptSym}[0]{\ensuremath{\DualPESym^{\C*}}}

\newcommand{\DualPECLbSym}[0]{\ensuremath{\lbsym{\DualPESym}^{\C}}}
\newcommand{\DualPECLbOpt}[2]{\ensuremath{\lbsym{\DualPESym}^{\C*}_{({#1},{#2})}}}
\newcommand{\DualPECUbSym}[0]{\ensuremath{\ubsym{\DualPESym}^{\C}}}

\newcommand{\DualPECUbOpt}[2]{\ensuremath{\ubsym{\DualPESym}^{\C*}_{({#1},{#2})}}}

\newcommand{\DualOverSym}[0]{\ensuremath{\sigma}}

\newcommand{\DualOverCOptSym}[0]{\ensuremath{\DualOverSym^{\C*}}}

\newcommand{\DualOverCOpt}[2]{\ensuremath{\DualOverCOptSym_{({#1},{#2})}}}
\newcommand{\DualOverCLbSym}[0]{\ensuremath{\lbsym{\DualOverSym}^{\C}}}
\newcommand{\DualOverCLb}[2]{\ensuremath{\DualOverCLbSym_{({#1},{#2})}}}
\newcommand{\DualOverCLbOptSym}[0]{\ensuremath{\lbsym{\DualOverSym}^{\C*}}}
\newcommand{\DualOverCLbOpt}[2]{\ensuremath{\DualOverCLbOptSym_{({#1},{#2})}}}
\newcommand{\DualOverCUbSym}[0]{\ensuremath{\ubsym{\DualOverSym}^{\C}}}
\newcommand{\DualOverCUb}[2]{\ensuremath{\DualOverCUbSym_{({#1},{#2})}}}
\newcommand{\DualOverCUbOptSym}[0]{\ensuremath{\ubsym{\DualOverSym}^{\C*}}}
\newcommand{\DualOverCUbOpt}[2]{\ensuremath{\DualOverCUbOptSym_{({#1},{#2})}}}

\newcommand{\DualSOnSym}[0]{\ensuremath{\alpha}}

\newcommand{\DualSOnCLbSym}[0]{\ensuremath{\lbsym{\DualSOnSym}^{\C}}}
\newcommand{\DualSOnCUbSym}[0]{\ensuremath{\ubsym{\DualSOnSym}^{\C}}}

\newcommand{\DualSOnCLbOptSym}[0]{\ensuremath{\lbsym{\DualSOnSym}^{\C*}}}
\newcommand{\DualSOnCUbOptSym}[0]{\ensuremath{\ubsym{\DualSOnSym}^{\C*}}}
\newcommand{\DualSOnCLbOpt}[2]{\ensuremath{\DualSOnCLbOptSym_{({#1},{#2})}}}
\newcommand{\DualSOnCUbOpt}[2]{\ensuremath{\DualSOnCUbOptSym_{({#1},{#2})}}}

\newcommand{\DualSOffSym}[0]{\ensuremath{\tau}}

\newcommand{\DualSOffCLbOptSym}[0]{\ensuremath{\lbsym{\DualSOffSym}^{\C*}}}
\newcommand{\DualSOffCLbOpt}[2]{\ensuremath{\DualSOffCLbOptSym_{({#1},{#2})}}}
\newcommand{\DualSOffCUbOptSym}[0]{\ensuremath{\ubsym{\DualSOffSym}^{\C*}}}
\newcommand{\DualSOffCUbOpt}[2]{\ensuremath{\DualSOffCUbOptSym_{({#1},{#2})}}}
\newcommand{\DualSOffCLbSym}{\ensuremath{\lbsym{\DualSOffSym}^{\C}}}
\newcommand{\DualSOffCUbSym}{\ensuremath{\ubsym{\DualSOffSym}^{\C}}}

\newcommand{\DualVadOffSym}[0]{\ensuremath{\xi}}

\newcommand{\DualVadOffCLbSym}[0]{\ensuremath{\lbsym{\DualVadOffSym}^{\C}}}
\newcommand{\DualVadOffCUbSym}[0]{\ensuremath{\ubsym{\DualVadOffSym}^{\C}}}

\newcommand{\DualVadOffCLbOpt}[2]{\ensuremath{\lbsym{\DualVadOffSym}^{\C*}_{({#1},{#2})}}}
\newcommand{\DualVadOffCUbOpt}[2]{\ensuremath{\ubsym{\DualVadOffSym}^{\C*}_{({#1},{#2})}}}

\newcommand{\masterSol}{\DemandFixedSym, \zFixedSym, \ShiftFixedSym}
\newcommand{\subP}{\ensuremath{S_{\C}(\masterSol, \etaFixedSym_\C)}}

\newcommand{\lesP}{\ensuremath{L^{\masterSol, \C, p}}}

\newcommand{\ie}{i.e.}
\newcommand{\eg}{e.g.}

\newcommand{\wolog}{w.l.o.g.}

\definecolor{rwth-blue}{RGB}{0,84,159}
\definecolor{rwth-lblue}{RGB}{142,186,229}
\definecolor{rwth-purple}{RGB}{122,111,172}
\definecolor{rwth-violet}{RGB}{97,33,88}
\definecolor{rwth-carmine}{RGB}{161,16,53}
\definecolor{rwth-red}{RGB}{204,7,30}
\definecolor{rwth-magenta}{RGB}{227,0,102}
\definecolor{rwth-orange}{RGB}{246,168,0}
\definecolor{rwth-yellow}{RGB}{255,237,0}
\definecolor{rwth-grass}{RGB}{189,171,39}
\definecolor{rwth-green}{RGB}{87,171,39}
\definecolor{rwth-cyan}{RGB}{0,152,161}
\definecolor{rwth-teal}{RGB}{0,97,101}
\definecolor{rwth150red}{HTML}{a11035}

\DeclareDocumentCommand\colorsOn{}{
  \definecolor{visual-red}{RGB}{255,0,0}
  \definecolor{visual-blue}{RGB}{0,0,255}
  \definecolor{visual-green}{RGB}{0,128,0}
  \definecolor{visual-purple}{RGB}{220,0,220}
  \definecolor{visual-orange}{RGB}{255,128,0}
  \definecolor{visual-cyan}{RGB}{0,200,200}
}

\colorsOn

\usepackage[ruled,vlined]{algorithm2e}
\SetKwInput{KwData}{Input}
\SetKwInput{KwResult}{Output}
\SetKw{KwGoto}{goto}
\SetKw{KwWith}{with}
\SetKw{AND}{and}
\SetKw{OR}{or}
\SetKw{GOTO}{go to}
\SetKw{STOP}{STOP.}

\newcommand{\autocref}[1]{Constraint \eqref{#1}}

\definecolor{ibm-blue}{RGB}{100,143,255}
\definecolor{ibm-purple}{RGB}{120,94,240}
\definecolor{ibm-red}{RGB}{220,38,127}
\definecolor{ibm-orange}{RGB}{254,97,0}
\definecolor{ibm-yellow}{RGB}{255,176,0}

\colorlet{color1}{ibm-blue}
\colorlet{color2}{ibm-orange}
\colorlet{color3}{ibm-red}
\colorlet{color4}{ibm-purple}
\colorlet{color5}{ibm-yellow}

\newcommand{\Nodes}{\ensuremath{V}}
\newcommand{\Timesteps}{\ensuremath{T}}

\newcommand{\Transformers}{\ensuremath{\mathcal{T}}}
\newcommand{\TransGroups}{\ensuremath{\Transformers^\text{sym}}}

\newcommand{\Busbars}{\ensuremath{\Nodes^{w}}}              % busbar nodes
\newcommand{\Breakers}{\ensuremath{\Nodes^{f}}}             % breaker nodes
\newcommand{\BusbarsIn}[1]{\ensuremath{\Busbars({#1})}}     % busbars in substation
\newcommand{\BreakersIn}[1]{\ensuremath{\Breakers({#1})}}   % breakers in substation
\newcommand{\Couplers}{\ensuremath{\Arcs^{w}}}              % coupler arcs

\newcommand{\GenNodesIn}[1]{\ensuremath{\Nodes^{\DemandParSym}({#1})}}

\newcommand{\ParArcs}{\ensuremath{\Arcs^\text{sym}}}

\newcommand{\MustRunGroups}{\ensuremath{\mathcal{G}^\text{mr}}}

\newcommand{\RedispatchGroups}{\ensuremath{\mathcal{G}^\text{rd}}}
\newcommand{\RedispatchShare}[3]{\ensuremath{\alpha_{{#1},{#2}}\tsup{#3}}}

\newcommand{\TapSym}{\ensuremath{\kappa}}

\newcommand{\TapMin}[2]{\ensuremath{\ubar{\TapSym}_{({#1},{#2})}}}
\newcommand{\TapMax}[2]{\ensuremath{\bar{\TapSym}_{({#1},{#2})}}}

\newcommand{\ShiftSym}{\ensuremath{\psi}}
\newcommand{\ShiftFixedSym}{\ensuremath{\fixedSym{\ShiftSym}}}
\newcommand{\Shift}[2]{\ensuremath{\ShiftSym_{({#1},{#2})}}}
\newcommand{\ShiftFixed}[2]{\ensuremath{\ShiftFixedSym_{({#1},{#2})}}}
\newcommand{\ShiftSlope}[2]{\ensuremath{R_{({#1},{#2})}}}
\newcommand{\ShiftConst}[2]{\ensuremath{S_{({#1},{#2})}}}

\newcommand{\MinSwitchTime}{\ensuremath{\ubar{T}}}

\newcommand{\MaxDefDev}{\ensuremath{\bar{n}^\text{dev}}}

\newcommand{\RampUp}[1]{\ensuremath{\DemandParSym^{+,\text{ramp}}_{#1}}}
\newcommand{\RampDown}[1]{\ensuremath{\DemandParSym^{-,\text{ramp}}_{#1}}}

\newcommand{\MinUpTime}[1]{\ensuremath{\ubar{T}^{\text{on}}_{#1}}}
\newcommand{\MinDownTime}[1]{\ensuremath{\ubar{T}^{\text{off}}_{#1}}}

\newcommand{\MaxRdCount}{\ensuremath{\bar{n}^{\text{rd}}}}

\newcommand{\DemandChangeIt}[2]{\ensuremath{\hat{\DemandSym}_{#1}\tsup{#2}}}
\newcommand{\TapChange}[3]{\ensuremath{\zSym^{\mathrm{tap}\tin{#3}}_{(#1,#2)}}}
\newcommand{\TapStepMax}[2]{\ensuremath{\Delta\TapSym^{\max}_{(#1,#2)}}}

\newcommand{\FlowRelDifMax}[2]{\ensuremath{\overline{\Delta L}_{({#1},{#2})}}}

\newcommand{\FlexThresh}{\ensuremath{\hat{\imath}^{\text{flex}}}}

\newcommand{\PotDifMaxCplAt}[1]{\ensuremath{\bar{\PotSym}^{\text{cp}}_{#1}}}

\newcommand{\CostRdP}[1]{\ensuremath{c^{+}_{#1}}}
\newcommand{\CostRdM}[1]{\ensuremath{c^{-}_{#1}}}
\newcommand{\CostCh}[2]{\ensuremath{c^{\text{ch}}_{({#1},{#2})}}}
\newcommand{\CostCp}[2]{\ensuremath{c^{\text{cp}}_{({#1},{#2})}}}
\newcommand{\CostSp}[2]{\ensuremath{c^{\text{sp}}_{({#1},{#2})}}}

\newcommand{\B}{\ensuremath{\{0,1\}}}

\newif\ifModelTime
\ModelTimetrue
\newcommand{\tsup}[1]{\ifModelTime^{#1}\fi}                  % standalone
\newcommand{\tin}[1]{\ifModelTime,#1\fi}                     % inside existing superscript
\newcommand{\qtime}{\ifModelTime,\;t\in\Timesteps\fi}        % \forall
\newcommand{\sumtime}{\ifModelTime\sum_{t\in\Timesteps}\fi}  % time sum

\newcommand{\zt}[3]{\ensuremath{\zSym_{({#1},{#2})}\tsup{#3}}}          % z_{(i,j)}^{t}
\newcommand{\zGent}[2]{\ensuremath{\zSym_{{#1}}\tsup{#2}}}              % z_{v}^{t}
\newcommand{\zChanget}[3]{\ensuremath{\zSym^{\text{ch}\tin{#3}}_{({#1},{#2})}}}
\newcommand{\zCouplet}[3]{\ensuremath{\zSym^{\text{cp}\tin{#3}}_{({#1},{#2})}}}
\newcommand{\zSplitt}[3]{\ensuremath{\zSym^{\text{sp}\tin{#3}}_{({#1},{#2})}}}
\newcommand{\zUpt}[2]{\ensuremath{\zSym^{\text{up}\tin{#2}}_{{#1}}}}
\newcommand{\zDownt}[2]{\ensuremath{\zSym^{\text{down}\tin{#2}}_{{#1}}}}

\newcommand{\Demandt}[2]{\ensuremath{\DemandSym_{{#1}}\tsup{#2}}}       % g_{v}^{t}
\newcommand{\DemandPlust}[2]{\ensuremath{\DemandSym^{+\tin{#2}}_{{#1}}}}
\newcommand{\DemandMinust}[2]{\ensuremath{\DemandSym^{-\tin{#2}}_{{#1}}}}

\newcommand{\Demandflext}[2]{\ensuremath{\DemandSym^{\text{flex}\tin{#2}}_{{#1}}}}

\newcommand{\Tapt}[3]{\ensuremath{\TapSym_{({#1},{#2})}\tsup{#3}}}
\newcommand{\Shiftt}[3]{\ensuremath{\ShiftSym_{({#1},{#2})}\tsup{#3}}}

\newcommand{\FlowCt}[3]{\ensuremath{\FlowSym^{\C\tin{#3}}_{({#1},{#2})}}}
\newcommand{\FlowOverCt}[3]{\ensuremath{\FlowOverSym_{({#1},{#2})}^{\C\tin{#3}}}}
\newcommand{\PotCt}[2]{\ensuremath{\PotSym^{\C\tin{#2}}_{{#1}}}}

\newcommand{\FlowRelDift}[3]{\ensuremath{\Delta L_{({#1},{#2})}\tsup{#3}}}

\newcommand{\Flowflext}[3]{\ensuremath{\FlowSym^{\text{flex}\tin{#3}}_{({#1},{#2})}}}
\newcommand{\Potflext}[2]{\ensuremath{\PotSym^{\text{flex}\tin{#2}}_{{#1}}}}

\newcommand{\GenDistt}[2]{\ensuremath{d^{\text{dist}\tin{#2}}_{{#1}}}}

\newcommand{\StatDt}[2]{\ensuremath{\StatDSym_{{#1}}\tsup{#2}}}
\newcommand{\DemandCostst}[2]{\ensuremath{\DemandCostsSym_{{#1}}\tsup{#2}}}
\newcommand{\DemandMint}[2]{\ensuremath{\ubar{\DemandParSym}\tsup{#2}_{{#1}}}}
\newcommand{\DemandMaxt}[2]{\ensuremath{\bar{\DemandParSym}\tsup{#2}_{{#1}}}}
\newcommand{\DemandZerot}[2]{\ensuremath{\DemandParSym^{0\tin{#2}}_{{#1}}}}
\newcommand{\DemandMinRdt}[2]{\ensuremath{\bar{\DemandParSym}^{-\tin{#2}}_{{#1}}}}
\newcommand{\DemandMaxRdt}[2]{\ensuremath{\bar{\DemandParSym}^{+\tin{#2}}_{{#1}}}}

\newcommand{\MustRunMint}[1]{\ensuremath{\ubar{\DemandParSym}^{\text{mr}\tin{#1}}}}
\newcommand{\MustRunMaxt}[1]{\ensuremath{\bar{\DemandParSym}^{\text{mr}\tin{#1}}}}
\newcommand{\DefDevCountt}[1]{\ensuremath{n^{\text{dev}\tin{#1}}}}

\newcommand{\ScoreSym}{\ensuremath{s}}
\newcommand{\WeightSym}{\ensuremath{w}}
\newcommand{\Score}[3]{\ensuremath{\ScoreSym^{\text{#1}}_{\text{#2}}({#3})}}   % s^{sec}_{name}(elem)
\newcommand{\ScoreTot}[2]{\ensuremath{\ScoreSym^{\text{#1}}({#2})}}           % s^{sec}(elem)
\newcommand{\Weight}[2]{\ensuremath{\WeightSym^{\text{#1}}_{\text{#2}}}}       % w^{sec}_{name}

\newcommand{\Viol}[2]{\ensuremath{\nu({#1},{#2})}}
\newcommand{\ObjMP}[1]{\ensuremath{\Phi_{\text{mp}}({#1})}} 
\newcommand{\ObjInc}[1]{\ensuremath{\Phi_{\text{inc}}({#1})}}
\newcommand{\LODF}[2]{\ensuremath{\Lambda_{{#1},{#2}}}}
\newcommand{\HopDist}[2]{\ensuremath{d({#1},{#2})}}

\newcommand{\spmaxn}{\ensuremath{\bar{k}^{\text{sp}}}}
\newcommand{\spmaxrel}{\ensuremath{\bar{\alpha}^{\text{sp}}}}
\newcommand{\spminscore}{\ensuremath{\ubar{\ScoreSym}^{\text{sp}}}}
\newcommand{\sppb}{\ensuremath{\beta^{\text{sp}}}}
\newcommand{\spminfeas}{\ensuremath{\ubar{k}^{\text{sp}}_{\text{feas}}}}
\newcommand{\spminopt}{\ensuremath{\ubar{k}^{\text{sp}}_{\text{opt}}}}
\newcommand{\spstab}{\ensuremath{k^{\text{sp}}_{\text{stab}}}}
\newcommand{\spdecay}{\ensuremath{\gamma^{\text{sp}}}}
\newcommand{\spsimk}{\ensuremath{k^{\text{sp}}_{\text{div}}}}

\newcommand{\spItSolved}[2]{\ensuremath{I_{{#1}{#2}}}}
\newcommand{\spItCut}[2]{\ensuremath{C_{{#1}{#2}}}}

\newcommand{\spCuts}[1]{\ensuremath{K_{#1}}}
\newcommand{\spCutsFeas}[1]{\ensuremath{K_{{#1},\text{feas}}}}
\newcommand{\spCutsOpt}[1]{\ensuremath{K_{{#1},\text{opt}}}}

\newcommand{\spS}[2]{\ScoreTot{sp}{{#1},{#2}}}
\newcommand{\spSrel}[2]{\Score{sp}{rel}{{#1},{#2}}}
\newcommand{\spSshr}[2]{\Score{sp}{shr}{{#1},{#2}}}
\newcommand{\spSvio}[2]{\Score{sp}{vio}{{#1},{#2}}}
\newcommand{\spSstl}[2]{\Score{sp}{stl}{{#1},{#2}}}
\newcommand{\spSrnd}[2]{\Score{sp}{rnd}{{#1},{#2}}}
\newcommand{\spSdiv}[2]{\Score{sp}{div}{{#1},{#2}}}
\newcommand{\spWrel}{\Weight{sp}{rel}}
\newcommand{\spWshr}{\Weight{sp}{shr}}
\newcommand{\spWvio}{\Weight{sp}{vio}}
\newcommand{\spWstl}{\Weight{sp}{stl}}
\newcommand{\spWrnd}{\Weight{sp}{rnd}}
\newcommand{\spWdiv}{\Weight{sp}{div}}

\newcommand{\cfmaxn}{\ensuremath{\bar{k}^{\text{cut}}}}
\newcommand{\cfmaxrel}{\ensuremath{\bar{\alpha}^{\text{cut}}}}
\newcommand{\cfpb}{\ensuremath{\beta^{\text{cut}}}}
\newcommand{\cfn}[1]{\ensuremath{k^{\text{cut}}({#1})}}
\newcommand{\cfSvio}[2]{\Score{cut}{vio}{{#1},{#2}}}

\newcommand{\pdcphop}{\ensuremath{\bar{d}^{\text{pd}}_{\text{cp}}}}

\newcommand{\pdmaxn}{\ensuremath{\bar{k}^{\text{pd}}}}
\newcommand{\pdmaxrel}{\ensuremath{\bar{\alpha}^{\text{pd}}}}
\newcommand{\pdflow}[1]{\ensuremath{\fixedSym{\FlowSym}^{\text{pd}}_{#1}}}
\newcommand{\pdcong}[0]{\ensuremath{\Arcs^{\text{pd}}_{\text{cong}}}}
\newcommand{\pdinc}{\ensuremath{\aC^{\text{pd}}_{\text{inc}}}}
\newcommand{\pddec}{\ensuremath{\aC^{\text{pd}}_{\text{dec}}}}

\newcommand{\pdS}[1]{\ScoreTot{pd}{#1}}
\newcommand{\pdScp}[1]{\Score{pd}{cp}{#1}}
\newcommand{\pdSsns}[1]{\Score{pd}{sns}{#1}}
\newcommand{\pdSvio}[1]{\Score{pd}{vio}{#1}}

\newcommand{\pdWcp}{\Weight{pd}{cp}}
\newcommand{\pdWsns}{\Weight{pd}{sns}}
\newcommand{\pdWvio}{\Weight{pd}{vio}}

\newcommand{\wsmaxn}{\ensuremath{\bar{k}^{\text{sh}}}}
\newcommand{\wsmaxrel}{\ensuremath{\bar{\alpha}^{\text{sh}}}}
\newcommand{\wsminscore}{\ensuremath{\ubar{\ScoreSym}^{\text{sh}}}}
\newcommand{\wsarcs}{\ensuremath{\Arcs^{\text{sh}}_{\text{sw}}}}

\newcommand{\wsSadm}[1]{\Score{sh}{adm}{#1}}
\newcommand{\wsWcp}{\Weight{sh}{cp}}
\newcommand{\wsWsns}{\Weight{sh}{sns}}
\newcommand{\wsWvio}{\Weight{sh}{vio}}

\newcommand{\cmark}{\ding{51}}%
\newlength{\normalparindent}
\AtBeginDocument{\setlength{\normalparindent}{\parindent}}

\AtBeginDocument{\let\cite\citep}

\usepackage{geometry}
\usepackage{authblk}

\renewcommand{\orcidicon}[1]{\href{https://orcid.org/#1}{%
  \begin{tikzpicture}[xscale=0.025,yscale=-0.025,transform shape]\pic{orcidlogo};\end{tikzpicture}}}
\let\textcite\relax

\usepackage[backend=biber,style=apa,natbib=true]{biblatex}
\newenvironment{paperfrontmatter}{%
  \title{\papertitle}
  \author{Tim Donkiewicz\,\protect\orcidicon{0000-0002-5721-3563}}
  \author{Oliver Gaul\,\protect\orcidicon{0000-0002-2131-1911}}
  \author{Marco L\"ubbecke\,\protect\orcidicon{0000-0002-2635-0522}}
  \affil{Chair of Operations Research, RWTH Aachen University\\
    Kackertstra\ss e 7, 52072 Aachen, Germany}
  \date{}
  \maketitle
}{%
  \par\smallskip
  \begingroup
    \small
    \noindent\textbf{Equal contribution.} Tim Donkiewicz and Oliver Gaul contributed equally to this work.
    \par\smallskip\noindent\textbf{Corresponding author.} Tim Donkiewicz: \href{mailto:donkiewicz@or.rwth-aachen.de}{donkiewicz@or.rwth-aachen.de}.
    \par\smallskip\noindent\textbf{Funding.} \paperfunding
    \par
  \endgroup
}
\newenvironment{paperkeywords}{\par\smallskip\noindent\textbf{Keywords:} }{\par\medskip}
\newcommand{\sep}{\unskip; }

\newcommand{\papertitle}{Integrated Topology and Redispatch Optimization of Real-World Power Transmission Systems using Branch-and-Benders-Cut}
\newcommand{\paperfunding}{This work was supported by the German Federal Ministry for Economic Affairs and Climate Action (BMWK) under grants 03EI4056C and 03EI4056D within the project ``InnOpTEM - Innovative Ansätze zur Optimierung von Topologiemaßnahmen im Engpassmanagement des Netzbetriebs'' (September 2022--August 2025).}

\hypersetup{pdftitle={Integrated Topology and Redispatch Optimization of Real-World Power Transmission Systems using Branch-and-Benders-Cut},pdfauthor={Tim Donkiewicz, Oliver Gaul, Marco Luebbecke},pdfkeywords={OR in energy, transmission topology optimization, redispatch, Benders decomposition, computational study}}
\begin{document}
\begin{paperfrontmatter}
\begin{abstract}
Transmission system operators can reduce congestion management costs by jointly optimizing network topology, generator redispatch, and transformer settings.
The large number of switching decisions, detailed operational constraints, and security requirements under outages make this problem computationally demanding.
We present a full mixed-integer linear programming formulation that combines line switching, substation reconfiguration, generator commitment and redispatch, and transformer settings with operational requirements.
We solve the problem using an exact branch-and-Benders-cut algorithm.
We integrate specialized subproblem solution methods, subproblem selection, cut filtering, partial decomposition, and techniques for enforcing network connectivity.
A sequential switching heuristic provides initial solutions for the exact solve.
An extensive computational study on 130 synthetic instances demonstrates substantial reductions in solution time compared with standard Benders decomposition and shows how heuristic warmstarting improves solution quality, especially when more switching actions are allowed.
Results on real transmission networks illustrate how joint topology and redispatch decisions identify effective switching actions that sequential approaches can miss.
Selected solutions generated by our approach were presented to system operators and subsequently checked using alternating current power flow calculations.
All seven discussed solutions passed a first operator assessment, indicating the practical potential of the integrated approach.
\end{abstract}
\begin{paperkeywords}
OR in energy \sep transmission topology optimization \sep redispatch \sep Benders decomposition \sep computational study
\end{paperkeywords}
\end{paperfrontmatter}

\section*{Introduction}\label{sec:introduction}
High-voltage power transmission grids enable the transfer of power between large scale power plants, consumers and lower-voltage network layers.
Here, \emph{congestion} refers to a (real or predicted) operation of individual power lines beyond their capacity, and the resulting need for, potentially costly countermeasures.
Annual congestion management costs range in the billions of euros~\cite{bnetza2024monitoring}, caused by,  among other things, a growing share of renewable, weather-dependent generation~\cite{titz2024identifying}.

Transmission system operators (TSOs) address these congestions in a multi-stage planning process.
One of these stages is the preventive redispatch process~(pRD2), which runs on a day-ahead basis.
During this step, a schedule of measures for the next day is planned.
This includes preventive measures for potential congestion in scenarios representing the failure of individual elements (N\nobreakdash-1-security).
The measures available during the congestion management process can be classified either as \emph{network-} or \emph{market-related}.
Market-related measures involve the redispatch, where pre-planned generator schedules are altered, which can incur high costs.
Network-related measures cover, for example, switching the network's topology or phase-shifting transformers, are relatively inexpensive, and have to be prioritized under German law~(\citealp[\S~13(1)]{enwgSystemResponsibility}; \citealp[p.~57]{bundestag2004EnwgDraft}).
However, only generator redispatch and phase-shifting transformer taps are optimized algorithmically \cite{ewerszumrode2024iterative}.
Topology measures, despite their demonstrated potential for congestion relief, are often handled manually and thus limited in both scope and quality.
Where they are determined algorithmically, they are often not optimized in conjunction with re-dispatch measures, which concedes additional improvement potential.
Integrating topology switching into an automated optimization procedure therefore promises significant benefits for grid security and cost efficiency.

Developing such optimization methods is challenging because of the large size of typical instances and the numerous operational and safety-related constraints that must be modelled and respected within topology optimization.
In this paper we report on our progress in this field.
It results from the InnOpTEM research project, which studied the integrated optimization of redispatch and topology measures in congestion management.
The project involved partners from research (RWTH Aachen University, Gurobi, FGH) and the four German transmission system operators (Amprion, 50Hertz Transmission, TenneT~TSO, TransnetBW).
This constellation enabled us to identify practically relevant and sufficient constraints, and test solution methods against real-world instances and time requirements.
In this paper, we develop a mixed-integer linear program (MILP) that we exactly solve with a tailored branch-and-Benders-cut algorithm.
For those interested in the scope, complementary approaches, and results of the project, we refer to the final report~\cite{innoptem-abschlussbericht}.

\paragraph{Problem and solution approach}
\label{subsec:problem-approach}

The optimization problem is an extended version of the security-constrained optimal transmission switching problem (SCOTS) with N\nobreakdash-1 scenarios.
The decision variables comprise the state of switchable branches and busbars (substation reconfiguration), the tap position of phase-shifting transformers, and the active-power dispatch (or redispatch) of generators.
Power flows are modeled with a lossless DC approximation and must satisfy thermal limits and angle bounds in the base case as well as under each N\nobreakdash-1 scenario.

The compact MILP replicates the entire DC power-flow system for every scenario.
For a realistic instance representing a national grid, with several thousand buses and hundreds of N\nobreakdash-1 scenarios, this model can already exceed the available memory of standard workstation hardware, so a monolithic solve is not an option.
We therefore solve SCOTS via a branch-and-Benders-cut algorithm that places topology, redispatch, and switching variables in the master problem and assigns one subproblem to each contingency.
Additions to the algorithm tailored to the use case, as presented in the following, make the instances computationally tractable.

\paragraph{Contribution}
\label{subsec:contributions}

This paper presents an integrated optimization approach to a SCOTS problem that realistically captures many of the operational requirements of German TSOs.
Its contribution lies in making the full model explicit, outlining a practically successful solution method, and documenting the advantages of an exact MILP algorithm over heuristic approaches.
The methodology is also a novel contribution for more general SCOTS problems.
The specific contributions are:
\begin{enumerate}
	\item A full mathematical formulation of the SCOTS problem established within the research project InnOpTEM~\cite{innoptem-abschlussbericht} with the German TSOs, including substation reconfiguration, phase-angle difference limits for disabled branches, proportional redispatch groups and transformer groups with equal tap positions, switch-back limits, minimum on/off times, and differentiated base/contingency thermal ratings (\Cref{sec:model}). That section is split between the ``core'' problem and extensions which include additional, practically relevant constraints. In parts, it resembles common literature formulations. The contribution here lies mainly in the full collection of constraints rather than any individual set of constraints.
	\item A branch-and-Benders-cut (B\&BC) algorithm that integrates the subproblem solution and cut-filtering methods developed in companion papers~\citep{donkiewicz2026potflow,donkiewicz2026cutfilter}, combines cut filtering with subproblem selection for SCOTS, and includes additional components such as partial decomposition and a connectivity separation procedure, as described in \Cref{sec:methodology}. We also discuss a sequential switching heuristic over restricted switching decisions, which prioritizes the practically essential primal solution quality.  This section focuses on solving the problem exactly, and can be applied to other problems that are sufficiently similar to the core problem.
	\item An extensive joint computational study of the integrated algorithm on instances generated with the methodology of~\textcite{donkiewicz2026instances}. In \Cref{sec:evaluation}, we compare selected algorithm configurations in terms of solution time and primal solution quality. We also comment on our evaluation on real-world instances, though we are unable to share concrete instances or results.
\end{enumerate}

\paragraph{Paper outline}
\label{subsec:outline}

The remainder of the paper is organized as follows.
\Cref{sec:related-work} positions our work against existing literature solution approaches.
\Cref{sec:model} gives our full problem formulation, with \Cref{subsec:dc-ots} introducing the core problem, and \Cref{subsec:extensions} providing additional, practically relevant extensions.
The core problem in \Cref{subsec:dc-ots} suffices as background for the remaining paper.
In \Cref{sec:methodology}, we discuss our extended branch-and-Benders-cut algorithm.
Computational results are reported in \Cref{sec:evaluation}, and operational case studies, operator assessments, and practical implications in \Cref{sec:practice}.

\section{Related Work}\label{sec:related-work}
\subsection{Optimal Transmission Switching}
\label{subsec:ots-lit}

The optimal transmission switching (OTS) problem extends the optimal power flow (OPF) problem by not only choosing the generation dispatch, but also deciding which branches (i.e., the edges of the underlying graph) to switch on or off~\cite{fisher2008optimal}.
Different formulations are used to approximate the alternating current (AC) power flow.
Lossless direct-current (DC) approximations represent active power flows through voltage-angle differences, omitting reactive power and voltage-magnitude variation~\cite{stott2009dc}.
Security-constrained OTS additionally accounts for outages, such as the failure of a branch~\cite{hedman2009optimal}.
Switching has also been combined with security-constrained unit commitment, which introduces generator on/off decisions over a planning horizon~\cite{hedman2010cooptimization,khodaei2010transmission}, additionally introducing ramps and startup/shutdown constraints.

\emph{Substation reconfiguration} changes the connections within substations by assigning lines, generators, and loads to busbars or opening and closing connections between busbars~\cite{morsy2022security}.
\textcite{heidarifar2016network} model the individual breakers within substations, allowing both line switching and substation reconfiguration.
\textcite{heidarifar2021optimal} model line switching and substation reconfiguration, with AC feasibility checks for the base case and N\nobreakdash-1 contingencies.
\textcite{bastianel2026switching} model line switching and substation reconfiguration in networks with AC and DC transmission, and compare a nonlinear formulation with convex relaxations and linear approximations.
\textcite{rajaei2026security} consider outages of lines, busbars, and couplers, using a linearized AC power-flow model and a heuristic decomposition with separate master problems for dispatch and substation reconfiguration.
Operational constraints considered in the literature include phase-angle difference limits for disabled branches~\cite{khodaei2012voltage} and connectivity requirements, enforced through an auxiliary flow LP whose optimal value distinguishes connected from disconnected topologies~\cite{han2022connectedness}.

\subsection{Formulations and Decomposition Approaches}
\label{subsec:benders-lit}

Formulation strength and decomposition both affect the computational effort required to solve OTS.
For DC OTS, \citet{pineda2024tight} perform optimization-based big-$M$ tightening.
\textcite{kocuk2016cycle} develop an alternative cycle-based formulation.
\textcite{goldis2017shift} use shift factors, which express how changes in generation minus demand at each node affect line flows, to formulate security-constrained line switching and substation reconfiguration compactly.

Benders decomposition splits variables/decisions into a master problem (MP) and subproblems (SPs) linked by feasibility and optimality cuts~\cite{benders1962partitioning}.
Lazily separating the Benders cuts within the branch-and-bound tree yields a branch-and-Benders-cut algorithm~\cite{rahmaniani2017benders}.
Benders decomposition is used for security-constrained optimal transmission switching with unit commitment, where subproblems guarantee power-flow feasibility under outages~\cite{khodaei2010transmission,khodaei2012voltage,schumacher2017algorithm,sheikh2019security,ramesh2022accelerated}, and for joint line switching and substation reconfiguration~\cite{heidarifar2014reconfiguration}.
\textcite{morsy2022security} use column-and-constraint generation for security-constrained substation reconfiguration.

The structure of the power-flow model also permits specialized solution procedures.
For affine potential-based flow problems, \textcite{donkiewicz2026potflow} solve the Benders subproblem via linear equation systems and factorization updates.
Low-rank updates also accelerate the GPU-based DC power-flow calculations of \textcite{westerbeck2025accelerated} after topology changes.

With many contingencies, three further choices affect the solution process: which subproblems to solve (next), which generated cuts to add, and which contingency constraints to retain explicitly in the master.
\textcite{donkiewicz2026subproblem} studies adaptive subproblem selection for survivable network design, including scores trained online by logistic regression.
\textcite{donkiewicz2026cutfilter} study filtering of Benders cuts for affine potential-based flow problems.
Partial Benders decomposition, as studied by \textcite{crainic2021partial}, retains selected scenario variables and constraints in the master to strengthen it.

\subsection{Heuristic and Learning Approaches in Operational Settings}
\label{subsec:heuristic-lit}

Operational time limits motivate heuristics that reduce the search over topology decisions.
Early work includes sensitivity analysis~\cite{hedman2008sensitivity}, switching heuristics~\cite{fuller2012fast}, and heuristic prescreening of switchable branches~\cite{liu2012heuristic}.
\textcite{babaeinejadsarookolaee2023congestion} present a sensitivity-based heuristic for substation reconfiguration.

\textcite{lai2023optimisation} alternate unit commitment and topology optimization for DC security-constrained line switching and substation reconfiguration.
\textcite{ewerszumrode2024iterative} alternate topology optimization with redispatch and phase-shifter optimization using linearized models.
\textcite{preuschoff2024investigating} heuristically preselect topologies, then jointly optimize topology and redispatch over the retained candidates.

Learning methods can guide optimization or propose candidate actions.
\textcite{preuschoff2024machine} study machine learning heuristics for network topology optimization, and \textcite{pineda2024learning} use solutions of past DC OTS instances to guide new solves, without guaranteeing optimality.
\textcite{rajaei2025graph} use a graph neural network to select candidate substations for reconfiguration and evaluate transfer across network configurations and systems.

\subsection{Positioning of This Work}
\label{subsec:positioning}
We present a full mathematical formulation and an integrated solution approach for an operationally realistic topology optimization problem established with the German TSOs.
Our DC formulation combines line switching, substation reconfiguration, generator commitment and redispatch, and phase-shifting transformer settings with detailed operational constraints, including proportional redispatch groups, equal-tap transformer groups, and limits on switching changes (\Cref{sec:model}).
The same topology, generation, and transformer decisions apply to the base case and every contingency.
To our knowledge, this combination has not previously been published as a full MILP formulation at this level of operational detail.
\Cref{tab:related-work-comparison} compares selected approaches (including two additional approaches developed in the project) by model scope and solution method.

\begin{table}[tbp]
  \centering
  \caption{Selected approaches to transmission topology optimization. Ticks indicate modeled decisions and constraints; dashes indicate components outside the stated scope. Status refers to the complete solution procedure, including candidate selection.}
  \label{tab:related-work-comparison}
  \small
  \setlength{\tabcolsep}{3pt}
  \renewcommand{\arraystretch}{1.2}
  \begin{tabularx}{\linewidth}{@{}>{\raggedright\arraybackslash}p{0.30\linewidth} c c c c c c c >{\raggedright\arraybackslash}X l@{}}
    \toprule
    & Flow & D & UC & L & SR & P & N\nobreakdash-1 & Method & Status \\
    \midrule
    \tablecite{hedman2010cooptimization}
    & DC & \cmark & \cmark & \cmark & -- & -- & \cmark & Alternating MILPs & Heuristic \\
    \tablecite{kocuk2016cycle}
    & DC & \cmark & -- & \cmark & -- & -- & -- & Cycle cuts + MILP & Exact \\
    \tablecite{goldis2017shift}
    & DC & \cmark & -- & \cmark & \cmark$^{*}$ & -- & \cmark & Screening + MILP & Heuristic \\
    \tablecite{heidarifar2021optimal}
    & AC-a & \cmark & -- & \cmark & \cmark & -- & \cmark & Decomposition; AC checks & Heuristic \\
    \tablecite{lai2023optimisation}
    & DC & \cmark & \cmark & \cmark & \cmark & -- & \cmark & Alternating MILPs & Heuristic \\
    \tablecite{rajaei2026security}
    & AC-a & \cmark & -- & -- & \cmark & -- & \cmark & Multiple masters & Heuristic \\
    \addlinespace
    \tablecite{ewerszumrode2024iterative}
    & AC-a & \cmark & \cmark & \cmark & \cmark & \cmark & \cmark & Alternating MILPs & Heuristic$^{\ddagger}$ \\
    \tablecite{preuschoff2024investigating}
    & AC-a & \cmark & \cmark & \cmark & \cmark & \cmark & \cmark & Screening + MILP & Heuristic$^{\ddagger}$ \\
    \addlinespace
    \textbf{This work}
    & DC & \cmark & \cmark & \cmark & \cmark & \cmark & \cmark & Branch-and-Benders-cut & Exact \\
    \bottomrule
  \end{tabularx}
  \par\smallskip
  \begin{minipage}{\linewidth}
    \footnotesize
    D: dispatch; UC: generator on/off decisions; L: line switching; SR: substation reconfiguration; P: transformer phase shifting;
    N\nobreakdash-1: constraints for modeled single-element outages.
    AC-a: AC approximation or relaxation; Heidarifar et al. additionally check full AC feasibility.
    Exact means globally exact for the stated formulation, including DC approximations, at convergence; heuristic candidate restrictions do not qualify.
    $^{*}$Goldis et al. use specified substation configurations without independent element assignments.
    $^{\ddagger}$Project implementations of iterative optimization and topology pre-selection, respectively; scope follows~\textcite{innoptem-abschlussbericht}. Both restrict topology choices and contingency constraints; the iterative method can add newly critical cases. In pre-selection, commitment constraints are enforced after fixing the topologies.
  \end{minipage}
\end{table}

Our branch-and-Benders-cut algorithm integrates the subproblem algorithms and cut-filtering methods developed in the companion papers with subproblem selection for SCOTS, problem-specific partial decomposition, and supporting techniques (\Cref{sec:methodology}).
Compared with alternating topology and redispatch optimization~\cite{ewerszumrode2024iterative} and methods restricted to preselected topologies~\cite{preuschoff2024investigating}, our method jointly optimizes the decisions admitted by the implemented DC formulation within a branch-and-Benders-cut framework.
A sequential heuristic uses restricted switching problems to find feasible solutions quickly, and improve them further in the main optimization problem.

The MILP model we present is not only comprehensive, we also solve it exactly. Next to these two distinguishing aspects, we would like to mention that the computational literature in this area oftentimes experiments with only few instances~\citep{donkiewicz2026instances}. By contrast,
our computational study systematically tunes and integrates algorithmic components, evaluated on 130 synthetic test instances.
As a result, for a single allowed switching action, the integrated algorithm reduces the shifted geometric mean solve time by \qty[round-precision=1]{88.1}{\percent} relative to basic Benders decomposition with LP subproblems (\Cref{tab:final-parameter-study}).
With larger switching limits, it finds further objective reductions through solutions using multiple switching actions.

We also report feasible solutions on real transmission-network instances under the implemented operational constraints, reproducing evidence from the final project report~\cite{innoptem-abschlussbericht}.
Seven selected solutions, including line switching and substation reconfiguration, passed operator assessments and separate AC feasibility checks.
We relate these results to the practical benefits and limitations of integrated optimization (\Cref{sec:practice}).

\section{Model}\label{sec:model}
% Section 3: Model
% Security-Constrained Optimal Transmission Switching with Redispatch.
%
% Notation: reuses innoptem graph-theoretic macros from common-resources.tex,
%           extended by sets_and_defs.tex.

We formulate security-constrained optimal transmission switching (SCOTS) as a mixed-integer linear program over a DC power flow approximation~\cite{stott2009dc}.
The model decides on branch switching states, generator commitment and redispatch, transformer tap positions, and substation reconfiguration across a discrete consecutive planning horizon~$\Timesteps \subseteq \N$, subject to N\nobreakdash-1 contingency security.
\Cref{subsec:dc-ots} describes the core problem formulation, which is the most relevant part for the following sections, and \Cref{subsec:extensions} outlines additional, practically relevant extensions.
As stated earlier, some constraint groups here already have related formulations in (but are scattered over) the literature, which we indicate in the paragraph headers.

%% ======================================================================
\subsection{Base Model}
\label{subsec:dc-ots}
%% ======================================================================

Our base model is a single-time step DC-SCOTS formulation with transformer shifting, generator commitment, and substation reconfiguration.

We  represent the power network by a directed graph $G = (\Nodes,\Arcs)$ with \emph{buses} $\Nodes$ and \emph{arcs} $\Arcs \subseteq \Nodes \times \Nodes$.
Every node $u\in\Nodes$ carries a fixed \emph{load}~$\StatDt{u}{t}\in\R$ (positive for consumption) and a controllable \emph{generation}~$\Demandt{u}{t}\in\R$ (positive for infeed), bounded by $\DemandMint{u}{t} \in \R$ and $\DemandMaxt{u}{t} \in \R$.
Arcs are further grouped into \emph{branches} $\LArcs \subseteq \Arcs$, with detailed modeling of power flow, and \emph{substation arcs}, $\SArcs \subseteq \Arcs$, which, as elaborated later in this section, are virtual arcs that model substation topology.
Any arc is part of one or both of these sets, and has a variable \emph{switching state} $\zt{u}{v}{t}$, which might be fixed.
If all switching states are fixed, the problem is a \emph{security-constrained optimal power flow problem (SCOPF)}.

Symbols referring to the ``previous'' time step $t \leq 0$, $t \notin \Timesteps$ are parameters, \ie, they reflect the actions already taken in the previous time steps.

%% ======================================================================
\paragraph{\texorpdfstring{N\protect\nobreakdash-1}{N-1} Security~\textmd{\citep{hedman2009optimal}}}
\label{subsec:n1}
%% ======================================================================

The \emph{contingency set} $\aC \subseteq 2^{\LArcs}$ lists branch-failure scenarios that must be covered.
A branch that fails no longer influences the power flow on the remaining network.
Since distinct contingencies may have differing power flows, a separate set of power flow variables is introduced for each contingency, as detailed in the next paragraph.
Generation, switching decisions, and tap positions are shared across all contingencies.
The contingency set contains the base case $\emptyset$ where no branch failures occur.
Where applicable, the following holds for some fixed but arbitrary contingency $\C \in \aC$.

\paragraph{Power flow and thermal limits~\textmd{\citep{fisher2008optimal,hedman2009optimal}}}
For each contingency a separate power flow is calculated.
It is determined by the interplay of \emph{voltage angle} variables $\PotCt{v}{t} \in \R$ at each node $u \in \Nodes$, and \emph{flow} variables $\FlowCt{u}{v}{t} \in \R$ at each arc $(u,v) \in \Arcs$.

For failed branches $(u,v)\in\C$, we fix $\FlowCt{u}{v}{t}=0$ for all $t\in\Timesteps$.
\begin{subequations}\label{eq:dc-pf}
\begin{alignat}{3}
	%% power balance
	\sum_{(v,u)\in\delta^-(u)\setminus\C}\!\FlowCt{v}{u}{t}
	- \sum_{(u,v)\in\delta^+(u)\setminus\C}\!\FlowCt{u}{v}{t}
	+ \Demandt{u}{t} - \StatDt{u}{t}
	&= 0
	&&\quad \forall\, u\in\Nodes\qtime,
	\label{eq:power-balance}
	\\[4pt]
	%% branch flow (big-M on)
	\bigl((\PotCt{u}{t}-\PotCt{v}{t}-\Shiftt{u}{v}{t})
		+ (1-\zt{u}{v}{t})\PotDifMaxOff{u}{v}\bigr)\Sus{u}{v}
	&\ge \FlowCt{u}{v}{t}
	&&\quad \forall\,(u,v)\in\LArcs\setminus\C\qtime,
	\label{eq:flow-ub}
	\\
	\bigl((\PotCt{u}{t}-\PotCt{v}{t}-\Shiftt{u}{v}{t})
		- (1-\zt{u}{v}{t})\PotDifMaxOff{u}{v}\bigr)\Sus{u}{v}
	&\le \FlowCt{u}{v}{t}
	&&\quad \forall\,(u,v)\in\LArcs\setminus\C\qtime,
	\label{eq:flow-lb}
	\\[4pt]
	%% thermal limits with overload slack + switching
	-\bigl(\FlowMaxOver{u}{v}+\FlowOverCt{u}{v}{t}\bigr)
	\;\le\; \FlowCt{u}{v}{t}
	&\le \FlowMaxOver{u}{v}+\FlowOverCt{u}{v}{t}
	&&\quad \forall\,(u,v)\in\LArcs\qtime,
	\label{eq:thermal_over}
	\\[4pt]
	-\zt{u}{v}{t}\FlowMax{u}{v}
	\;\le\; \FlowCt{u}{v}{t}
	&\le \zt{u}{v}{t}\FlowMax{u}{v}
	&&\quad \forall\,(u,v)\in\Arcs\qtime,
	\label{eq:thermal_sw}
\end{alignat}
\end{subequations}

Constraint~\eqref{eq:power-balance} enforces nodal power balance.
For every branch $(u,v) \in \LArcs\setminus\C$, Kirchhoff's law is implemented by Constraints~\eqref{eq:flow-ub}--\eqref{eq:flow-lb}, which link flow to the voltage angle difference via the susceptance~$\Sus{u}{v} \in \R_{> 0}$ and phase shift (see below)~$\Shiftt{u}{v}{t} \in \R$.
When $\zt{u}{v}{t}=0$, the big-$M$ terms $\Sus{u}{v}\PotDifMaxOff{u}{v}$ deactivate the coupling.
\autocref{eq:thermal_over} couples the flow with an overload variable~$\FlowOverCt{u}{v}{t}$ when the soft bounds $\FlowMaxOver{u}{v} \in \R_{\geq 0}$ are violated that is penalized with cost $\OverCost \in \R$ per unit in the objective.
For any arc $(u,v) \in \Arcs$, Constraint~\eqref{eq:thermal_sw} limits branch loading to the hard thermal capacity~$\FlowMax{u}{v} \in \R_{\geq 0}$ when the branch is turned on, and else forces it to $0$.
For substation arcs $\SArcs$, the limits are chosen large enough that they cannot be violated when the branch is on.

\paragraph{Transformer tap shifting}
Phase shifting transformers allow for offsetting the voltage angle difference of connected lines, and thus for more direct control of power flows.
We model them as part of branches, \ie, each branch $(u,v) \in \LArcs$ has additional parameters and variables that correspond to the parameters of the transformer situated between its endpoints.
The voltage angle offset depends on the transformer's tap position $\Tapt{u}{v}{t} \in \Z$: we can reach any tap position between $\TapMin{u}{v}, \TapMax{u}{v} \in \Z$, and one tap step changes the angle by $\ShiftSlope{u}{v} \in \R$ from the initial angle $\ShiftConst{u}{v} \in \R$.
\begin{subequations}\label{eq:tap}
\begin{alignat}{3}
	\Shiftt{u}{v}{t}
	&= \Tapt{u}{v}{t}\cdot\ShiftSlope{u}{v} + \ShiftConst{u}{v}
	&&\quad \forall\,(u,v)\in\LArcs\qtime,
	\label{eq:tap-shift}
	\\
	\TapMin{u}{v}
	\;\le\; \Tapt{u}{v}{t}
	&\le \TapMax{u}{v}
	&&\quad \forall\,(u,v)\in\LArcs\qtime.
	\label{eq:tap-bounds}
\end{alignat}
\end{subequations}
For non-transformer branches, $\ShiftSlope{u}{v} = \ShiftConst{u}{v} = 0$.

\paragraph{Voltage angle differences}
The DC approximation assumes $\sin(\Pot{u}-\Pot{v})\approx\Pot{u}-\Pot{v}$, which requires small angle differences (typically $\le \qty{30}{\degree}$).
We therefore bound the angle difference on every active arc $(u,v) \in \Arcs$ by $\PotDifMax{u}{v} \in \R$:
\begin{alignat}{3}
	-(1- \zt{u}{v}{t}) \PotDifMaxOff{u}{v} - \PotDifMax{u}{v}
	\le
	\PotCt{u}{t}-\PotCt{v}{t}-\Shiftt{u}{v}{t}
	&\le
	(1 - \zt{u}{v}{t}) \PotDifMaxOff{u}{v} + \PotDifMax{u}{v}
	&&\quad \forall\,(u,v)\in\LArcs\setminus\C\qtime,
	\label{eq:angle-diff}
\end{alignat}
with $\PotDifMaxOff{u}{v} \in \R$ as given for Constraint \eqref{eq:angle-off}.
One node $v_{\text{ref}} \in \Nodes$ is designated as the reference node with $\PotCt{v_{\text{ref}}}{t}=0$ for every time step $t \in \Timesteps$ to fix the angle origin.

%% ======================================================================
\paragraph{Substation reconfiguration~\textmd{\citep{rajaei2026security,bastianel2026switching}}}
\label{subsec:substation}
%% ======================================================================
To facilitate a wider range of topologies, nodes are divided into groups, so-called \emph{substations}.
Every node is part of a single substation $\Sub \in \Subs$ out of the set of substations $\Subs \subseteq 2^V$.
It is either part of the internal \emph{busbar nodes} $\BusbarsIn{\Sub}$ of the substation, or of its outward facing \emph{breaker nodes} $\BreakersIn{\Sub}$.
The nodes of the substation form a component connected by substation arcs.
Busbar nodes are connected by \emph{coupler arcs} $\Couplers \subseteq \SArcs$, and breaker nodes are connected to busbar nodes by \emph{breaker arcs}, which are always directed from the breaker node to the busbar node.
Breaker nodes are also connected to external equipment such as branches, generators, and loads.

\begin{figure}[H]
	\centering
	\newsavebox{\substationlayoutbox}
	\savebox{\substationlayoutbox}{%
		\begin{tikzpicture}[scale=0.82, every node/.style={scale=0.82}]
			\draw (0,.075) -- (2,.075);   %e1
			\draw (0,-.575) -- (2,-.575); %e2

			%%v1
			\draw (2.5,.75) -- (2.5,-4.675); % line v1
			\draw[fill=lightgray] (2,- .1  ) -- (3,- .1  ) -- (3,  .25 ) -- (2,  .25 ) -- cycle; % e1
			\draw[fill=white    ] (2,- .4  ) -- (3,- .4  ) -- (3,- .75 ) -- (2,- .75 ) -- cycle; % e2
			\draw[fill=lightgray] (2,-2.525) -- (3,-2.525) -- (3,-2.875) -- (2,-2.875) -- cycle; % coupler bot
			\draw[fill=lightgray] (2,-3.175) -- (3,-3.175) -- (3,-3.525) -- (2,-3.525) -- cycle; % e3
			\draw[fill=white    ] (2,-3.825) -- (3,-3.825) -- (3,-4.175) -- (2,-4.175) -- cycle; % e4

			%%v2
			\draw (4.2,.75) -- (4.2,-4.675); % line v2
			\draw[fill=white    ] (3.7,- .1  ) -- (4.7,- .1  ) -- (4.7,  .25 ) -- (3.7,  .25 ) -- cycle; % e1
			\draw[fill=lightgray] (3.7,- .4  ) -- (4.7,- .4  ) -- (4.7,- .75 ) -- (3.7,- .75 ) -- cycle; % e2
			\draw[fill=lightgray] (3.7,-1.05 ) -- (4.7,-1.05 ) -- (4.7,-1.4  ) -- (3.7,-1.4  ) -- cycle; % coupler top
			\draw[fill=white    ] (3.7,-3.175) -- (4.7,-3.175) -- (4.7,-3.525) -- (3.7,-3.525) -- cycle; % e3
			\draw[fill=lightgray] (3.7,-3.825) -- (4.7,-3.825) -- (4.7,-4.175) -- (3.7,-4.175) -- cycle; % e4

			%% K
			\draw (4.7,-1.225) -- (5.1,-1.225) -- (5.1,-2.7) -- (3,-2.7);
			\draw[fill=white    ] (4.95,-1.55) -- (5.25,-1.55) -- (5.25,-2.35) -- (4.95,-2.35) -- cycle;

			%% connections B1 -- B2
			\draw (3, 0.075) -- (3.7, 0.075); % e1
			\draw (3,- .575) -- (3.7,- .575); % e2
			\draw (3,-3.35 ) -- (3.7,-3.35 ); % e3
			\draw (3,-4    ) -- (3.7,-4    ); % e4

			\draw (4.7,-3.35) -- (6.7,-3.35); %e3
			\draw (4.7,-4) -- (6.7,-4);       %e4

			%% text
			\draw (0.3,  0.3  ) node {$a_1$};
			\draw (0.3, -0.35 ) node {$a_2$};
			\draw (6.4, -3.125) node {$a_3$};
			\draw (6.4, -3.775) node {$a_4$};
			\draw (2.18,  .68 ) node {$u$};
			\draw (3.88,  .68 ) node {$v$};
			\draw (5.5, -1.95 ) node {$a^w_1$};
			\draw (1.88, -4.38  ) node {$a^f_{u,4}$};
			\draw (3.58, -4.38  ) node {$a^f_{v,4}$};
		\end{tikzpicture}%
	}
	% height shared by both subfigures
	\newlength{\substationpanelht}
	\setlength{\substationpanelht}{\dimexpr\ht\substationlayoutbox+\dp\substationlayoutbox\relax}
	\begin{subfigure}[t]{.48\textwidth}\centering
		\parbox[c][\substationpanelht][c]{\linewidth}{\centering\usebox{\substationlayoutbox}}
		\caption{Schematic substation layout.
				Gray boxes represent closed switches, white open ones.}
		\label{fig:substation-layout}
	\end{subfigure}\hspace*{.02\textwidth}%
	\begin{subfigure}[t]{.48\textwidth}\centering
		\tikzstyle{aux_node} = [diamond, draw]
		\tikzset{
			coupler_node/.style = {circle, draw, minimum size=#1},
			coupler_node/.default = 11
		}
		\parbox[c][\substationpanelht][c]{\linewidth}{\centering
			\begin{tikzpicture}[scale=0.82, transform shape]
				% alignment helpers
				\node (0,0)                  (center) {};
				\node [left =1.8 of center] (lcenter) {};
				\node [right=1.8 of center] (rcenter) {};

				% busbars (circles)
				\node [coupler_node, above=1.3 of center] (coup1) {};
				\node [coupler_node, below=1.3 of center] (coup2) {};

				% breaker buses (diamonds)
				\node [aux_node, above =.15 of lcenter] (aux1) {};
				\node [aux_node, above =.15 of rcenter] (aux2) {};
				\node [aux_node, below =.15 of lcenter] (aux3) {};
				\node [aux_node, below =.15 of rcenter] (aux4) {};

				% external arc stubs
				\node [left =1.4 of aux1.west] (L1_west) {};
				\node [right=1.4 of aux2.east] (L2_east) {};
				\node [left =1.4 of aux3.west] (L3_west) {};
				\node [right=1.4 of aux4.east] (L4_east) {};

				\draw [->] (L1_west) -- (aux1);
				\draw [<-] (aux2)    -- (L2_east);
				\draw [<-] (L3_west) -- (aux3);
				\draw [->] (aux4)    -- (L4_east);

				% substation branches: solid = on, dashed = off
				\draw [->]         (aux1)  -- (coup1);
				\draw [->, dashed] (aux3)  -- (coup1);
				\draw [->, dashed] (aux1)  -- (coup2);
				\draw [->]         (aux3)  -- (coup2);
				\draw [->]         (aux2)  -- (coup1);
				\draw [->, dashed] (aux4)  -- (coup1);
				\draw [->, dashed] (aux2)  -- (coup2);
				\draw [->]         (aux4)  -- (coup2);
				\draw [->, dashed] (coup1) -- (coup2);

				% labels
				\node [above right=.03 and 0.005 of L1_west.east] {$a_1$};
				\node [above left =.03 and 0.005 of L2_east.west] {$a_3$};
				\node [above right=.03 and 0.005 of L3_west.east] {$a_2$};
				\node [above left =.03 and 0.005 of L4_east.west] {$a_4$};
				\draw [left=.1 of center] node {$a^w_1$};
				\node [above=.01 of coup1.north] {$u$};
				\node [below=.01 of coup2.south] {$v$};
				\draw (0.8, 0.35  ) node {$a^f_{u,4}$};
				\draw (1.3, -1.27  ) node {$a^f_{v,4}$};
			\end{tikzpicture}%
		}
		\caption{Graph representation with breaker nodes
				(diamonds), busbar nodes (circles) and switched on and off arcs (solid/dashed).}
		\label{fig:substation-graph}
	\end{subfigure}
	\caption{Substation modeling. Busbars $u,v$, incoming branches $a_1,\dots,a_4$, coupler $a^w_1$, breaker switches $a^f$.}
	\label{fig:substation}
\end{figure}
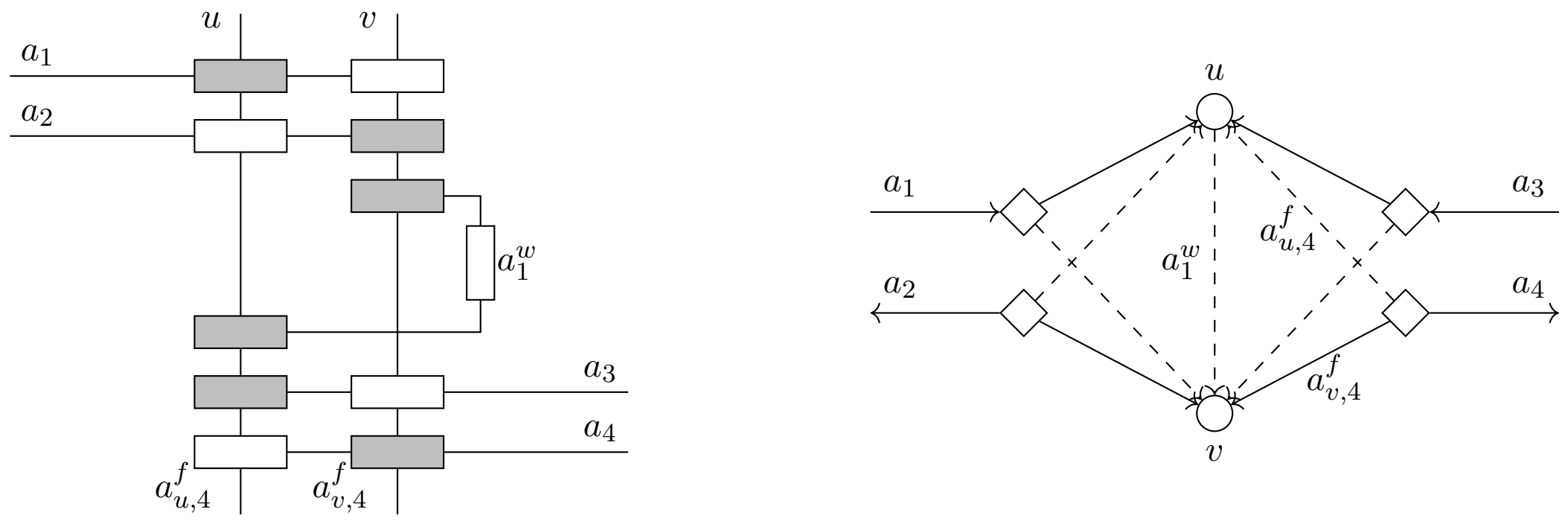

\Cref{fig:substation} illustrates a substation with busbars~$u$ and~$v$.
\Cref{fig:substation-layout} shows a classical schematic representation of a substation layout.
\Cref{fig:substation-graph} depicts the same layout in our graph representation.
Every incoming branch ends at a breaker node, which is then connected to every busbar by a substation-internal branch.
The coupler arc~$a^w_1$ directly links the two busbar nodes.

For the model, substation arcs are handled differently than branches.
We assume that the electrical distances covered are small enough that we can omit detailed modeling of the electrical parameters of the arcs.
As such, Kirchhoff's law \eqref{eq:flow-ub}--\eqref{eq:flow-lb} is not applied to substation arcs, but
simplified such that the voltage angles of incident buses need to be equal when the arc is switched on
\begin{subequations}
	\begin{equation}
		-(1 - \zt{u}{v}{t}) \PotDifMaxOff{u}{v} \leq \PotCt{u}{t} - \PotCt{v}{t} \leq (1 - \zt{u}{v}{t}) \PotDifMaxOff{u}{v} \quad \forall\,(u,v)\in\SArcs\setminus\C\qtime.
		\label{eq:substation-angle}
	\end{equation}
	Additionally, each breaker node may be connected to at most one busbar node
	\begin{alignat}{3}
		\sum_{v \in \BusbarsIn{\Sub}} \zt{u}{v}{t}
		&\le 1
		&&\quad \forall\,\Sub\in\Subs,\;
			u\in\BreakersIn{\Sub}\qtime.
		\label{eq:breaker}
	\end{alignat}
\end{subequations}

%% ======================================================================
\paragraph{Generator commitment~\textmd{\citep{hedman2010cooptimization}}}
\label{subsec:gen}
%% ======================================================================

At each node $v$ with a generator, the binary variable~$\zGent{v}{t}$ controls commitment.
Active power production is bounded by:
\begin{alignat}{3}
	\DemandMint{v}{t}\,\zGent{v}{t}
	\;\le\; \Demandt{v}{t}
	&\le \DemandMaxt{v}{t}\,\zGent{v}{t}
	&&\quad \forall\, v\in\Nodes\qtime.
	\label{eq:gen-bounds}
\end{alignat}

%% ======================================================================
\paragraph{Objective function}
\label{subsec:obj}
%% ======================================================================

The core objective minimizes production cost (coefficient $\DemandCostst{v}{t}$ per node) and overload penalty (weight $\OverCost$):
\begin{equation}\label{eq:objective}
	\min\;\;
	\underbrace{
	\sumtime
	\sum_{v\in\Nodes}
	\DemandCostst{v}{t}\,\Demandt{v}{t}
	}_{\text{production cost}}
	\;+\;
	\underbrace{
	\sumtime
	\sum_{\C \in \aC}
	\sum_{(u,v)\in\Arcs}
	\OverCost\,\FlowOverCt{u}{v}{t}
	}_{\text{overload penalty}}.
\end{equation}
In the extended model (\Cref{subsec:extensions}), the objective is augmented with redispatch and switching costs.

%% ======================================================================
\subsection{Model Extensions}
\label{subsec:extensions}
%% ======================================================================

The constraints above describe a single-time step SCOTS.
The following extensions add inter-temporal coupling, operational redispatch detail, and further practical requirements.

%% ------ Redispatch formulation ---------------------
\paragraph{Redispatch formulation}
\label{subsec:redispatch}
In the observed part of the planning process, generator dispatch has already been determined.
Due to branch capacities, it might however be necessary to redispatch generators, \ie, to deviate from the previous plan, which incurs additional costs, see Equation \eqref{eq:objective-ext}.
For each node $v \in \Nodes$ with planned generation setpoint $\DemandZerot{v}{t} \in \R$, variables $\DemandMinust{v}{t},\DemandPlust{v}{t} \in \R_{\geq 0}$ capture downward and upward deviations, bounded by $\DemandMinRdt{v}{t}, \DemandMaxRdt{v}{t} \in \R_{\geq 0}$.
A binary indicator $\DemandChangeIt{v}{t}\in\B$ records whether the generator is redispatched:
\begin{subequations}\label{eq:rd}
\begin{alignat}{3}
	\Demandt{v}{t}
	&= \DemandZerot{v}{t} + \DemandPlust{v}{t} - \DemandMinust{v}{t}
	&&\quad \forall\,v\in\Nodes,\;t\in\Timesteps,
	\label{eq:rd-link}
	\\
	0\;\le\;\DemandMinust{v}{t}
	&\le \DemandMinRdt{v}{t}\,\DemandChangeIt{v}{t}
	&&\quad \forall\,v\in\Nodes,\;t\in\Timesteps,
	\label{eq:rd-bounds}
	\\
	0\;\le\;\DemandPlust{v}{t}
	&\le \DemandMaxRdt{v}{t}\,\DemandChangeIt{v}{t}
	&&\quad \forall\,v\in\Nodes,\;t\in\Timesteps,
	\label{eq:rd-bounds-up}
	\\
	\sum_{t\in\Timesteps}\sum_{v\in\Nodes}\DemandChangeIt{v}{t}
	&\le \MaxRdCount.
	\label{eq:rd-count}
\end{alignat}
\end{subequations}
Constraint~\eqref{eq:rd-count} limits the number of generator--time-step pairs with redispatch over the planning horizon to~$\MaxRdCount\in\N$.
Note that, \wolog, we assume each node carries at most one generator.
In practice, multiple generating units at the same physical node can be modeled by introducing artificial auxiliary nodes (one per unit), each with its own cost and capacity parameters.

%% ------ Proportional redispatch groups ----------------
\paragraph{Proportional redispatch groups}
\label{subsec:redispatch-groups}
Some generators must be redispatched jointly in prescribed proportions.
Let $\RedispatchGroups\subseteq 2^\Nodes$ collect these groups, and let $\RedispatchShare{v}{G}{t} \in \R_{\ge 0}$ denote the fixed share of the generator at node $v\in \Nodes$ in group $G\in\RedispatchGroups \subseteq 2^\Nodes$, with $\sum_{v\in G}\RedispatchShare{v}{G}{t}=1$.
Each generator receives this share of the group's total upward and downward redispatch:
\begin{subequations}\label{eq:redispatch-group}
	\begin{alignat}{3}
		\DemandMinust{v}{t}
		&=\RedispatchShare{v}{G}{t}\sum_{w\in G}\DemandMinust{w}{t}
		&&\quad\forall\,G\in\RedispatchGroups,\;v\in G,\;t\in\Timesteps,
		\label{eq:redispatch-group-down}\\
		\DemandPlust{v}{t}
		&=\RedispatchShare{v}{G}{t}\sum_{w\in G}\DemandPlust{w}{t}
		&&\quad\forall\,G\in\RedispatchGroups,\;v\in G,\;t\in\Timesteps.
		\label{eq:redispatch-group-up}
	\end{alignat}
\end{subequations}
If these shares match the planned output proportions, all generators in the group change their output by the same percentage.
These constraints differ from the must-run group bounds~\eqref{eq:must-run}, which restrict aggregate output without fixing its distribution among generators.
When commitment decisions are enabled for a group, its members also share one on/off state:
\begin{equation}\label{eq:redispatch-group-status}
	\zGent{v}{t}=\zGent{w}{t}
	\quad\forall\,G\in\RedispatchGroups,\;v,w\in G,\;t\in\Timesteps.
\end{equation}
Their startup and shutdown indicators in~\eqref{eq:gen-updown} are shared as well.

%% ------ Busbar splitting and coupling -----------------
\paragraph{Busbar splitting and coupling~\textmd{\citep{lam2023multiperiod}}}
\label{subsec:busbar-coupling}
Variables $\zCouplet{u}{v}{t},\,\zSplitt{u}{v}{t}\in\B$ track whether coupler $(u,v)\in\Couplers$ is closed (``coupled'') or opened (``split'') between consecutive time steps $t-1$ and $t \in T$:
\begin{alignat}{3}
	\zt{u}{v}{t}
	&= \zt{u}{v}{t-1} - \zSplitt{u}{v}{t} + \zCouplet{u}{v}{t}
	&&\quad \forall\,(u,v)\in\Couplers,\;t\in\Timesteps.
	\label{eq:couple-link}
\end{alignat}

%% ------ Switching change tracking --------------------
\paragraph{Switching change tracking~\textmd{\citep{bastianel2025dayahead}}}
\label{subsec:switch-change}
For all arcs $(u,v) \in \Arcs$, the binary indicator $\zChanget{u}{v}{t}\in\B$ is~$1$ whenever the switching state changes between $t{-}1$ and~$t$:
\begin{subequations}\label{eq:change}
\begin{alignat}{3}
	\zt{u}{v}{t} - \zt{u}{v}{t-1}
	&\le \zChanget{u}{v}{t}
	&&\quad \forall\,(u,v)\in\Arcs,\;t\in\Timesteps,
	\label{eq:change-up}
	\\
	\zt{u}{v}{t-1} - \zt{u}{v}{t}
	&\le \zChanget{u}{v}{t}
	&&\quad \forall\,(u,v)\in\Arcs,\;t\in\Timesteps.
	\label{eq:change-down}
\end{alignat}
\end{subequations}

%% ------ Switching limits -----------------------------
\paragraph{Switching action limits~\textmd{\citep{fisher2008optimal,lam2023multiperiod,rajaei2026security}}}
\label{subsec:switch-limits}
The total number of switching actions is bounded by topology action limits $\TAL_t \in \Z$ (per time step), $\TAL_{(u,v)} \in \Z$ (per component), and $\TAL \in \Z$ (globally):
\begin{subequations}\label{eq:switch-lim}
\begin{alignat}{3}
	\sum_{(u,v)\in\Arcs}\zChanget{u}{v}{t}
	&\le \TAL_t
	&&\quad \forall\,t\in\Timesteps,
	\label{eq:switch-lim-t}
	\\
	\sum_{t\in\Timesteps}\zChanget{u}{v}{t}
	&\le \TAL_{(u,v)}
	&&\quad \forall\,(u,v)\in\Arcs,
	\label{eq:switch-lim-a}
	\\
	\sum_{t\in\Timesteps}\sum_{(u,v)\in\Arcs}\zChanget{u}{v}{t}
	&\le \TAL.
	\label{eq:switch-lim-all}
\end{alignat}
\end{subequations}

%% ------ Transformer adjustment limits ----------------
\paragraph{Transformer adjustment limits}
\label{subsec:tap-actions}
Transformer adjustments have separate action limits.
Let $\TapChange{u}{v}{t}\in\B$ indicate an adjustment of transformer $(u,v)$ in time step $t$, and let $\TapStepMax{u}{v}\geq 0$ bound the number of tap steps per adjustment:
\begin{equation}\label{eq:tap-change}
	-\TapStepMax{u}{v}\,\TapChange{u}{v}{t}
	\leq \Tapt{u}{v}{t}-\Tapt{u}{v}{t-1}
	\leq \TapStepMax{u}{v}\,\TapChange{u}{v}{t}
	\quad\forall\,(u,v)\in\Transformers,\;t\in\Timesteps.
\end{equation}
As for switching actions, limits apply per time step, per transformer, and over the full horizon:
\begin{subequations}\label{eq:tap-action-limits}
	\begin{alignat}{3}
		\sum_{(u,v)\in\Transformers}\TapChange{u}{v}{t}
		&\leq N^{\mathrm{tap}}_t
		&&\quad\forall\,t\in\Timesteps,\\
		\sum_{t\in\Timesteps}\TapChange{u}{v}{t}
		&\leq N^{\mathrm{tap}}_{(u,v)}
		&&\quad\forall\,(u,v)\in\Transformers,\\
		\sum_{t\in\Timesteps}\sum_{(u,v)\in\Transformers}\TapChange{u}{v}{t}
		&\leq N^{\mathrm{tap}}.
	\end{alignat}
\end{subequations}

%% ------ Switching continuity -------------------------
\paragraph{Minimum switching time~\textmd{\citep{lam2023multiperiod}}}
\label{subsec:min-switch-time}
To limit operational effort and equipment wear, a switching action cannot be reversed within $\MinSwitchTime \in \N$~time steps:
\begin{alignat}{3}
	\sum_{t'=t-\MinSwitchTime}^{t}\zChanget{u}{v}{t'}
	&\le 1
	&&\quad \forall\,(u,v)\in\Arcs,\;t\in\Timesteps.
	\label{eq:min-switch}
\end{alignat}

%% ------ Default topology -----------------------------
\paragraph{Default topology limit}
\label{subsec:default-topo}
Each arc $(u,v) \in \Arcs$ has a \emph{default} state~$\zDefault{u}{v} \in \B$.
The number of deviations is tracked by $\DefDevCountt{t}\in\Z_{\geq 0}$:
\begin{subequations}\label{eq:default}
\begin{alignat}{3}
	\bigl( \sum_{\substack{(u,v)\in\Arcs \\ \zDefault{u}{v}=1}}\zDefault{u}{v} - \zt{u}{v}{t} \bigr) + \bigl( \sum_{\substack{(u,v)\in\Arcs \\ \zDefault{u}{v}=0}} \zt{u}{v}{t} - \zDefault{u}{v} \bigr)
	&= \DefDevCountt{t}
	&&\quad \forall\,t\in\Timesteps,
	\label{eq:def-count}
	\\
	\DefDevCountt{t}
	&\le \MaxDefDev
	&&\quad \forall\,t\in\Timesteps.
	\label{eq:def-limit}
\end{alignat}
\end{subequations}
An upper deviation limit $\MaxDefDev \in \Z_{\geq 0}$ is imposed to prevent excessive deviation from the default topology, which may be desirable for operational familiarity and stability.

%% ------ Connectivity --------------------------------
\paragraph{Network connectivity~\textmd{\citep{han2020ensuring,han2022connectedness}}}
\label{subsec:connectivity}
The network has to be ``islanding-free.''
That is, for each contingency, all nodes with non-zero demands or generation not fixed to zero by their bounds have to be connected in the graph induced by the active arcs, \ie, those switched-on and not part of the contingency.
This is a weaker requirement than full connectivity, but ensures that all relevant nodes are always connected.

Let $\mathcal{K} \subset 2^\Arcs$ collect all minimal arc sets which, if switched off, violate this requirement in one of the contingencies.
For each $K\in\mathcal{K}$:
\begin{alignat}{3}
	\sum_{(u,v)\in K}\zt{u}{v}{t}
	&\ge 1
	&&\quad \forall\,t\in\Timesteps.
	\label{eq:connected}
\end{alignat}
$\mathcal{K}$ becomes intractably large for sufficiently large networks.
These constraints are added via a cut separation procedure, see \Cref{subsec:connectivity-method}.

%% ------ Generator startup / shutdown -----------------
\paragraph{Generator startup and shutdown~\textmd{\citep{hedman2010cooptimization}}}
\label{subsec:gen-startup}
Variables $\zUpt{v}{t},\,\zDownt{v}{t}\in\B$ track whether generator $v \in \Nodes$ is started up or shut down:
\begin{subequations}\label{eq:gen-updown}
\begin{alignat}{3}
	\zGent{v}{t}
	&= \zGent{v}{t-1} - \zDownt{v}{t} + \zUpt{v}{t}
	&&\quad \forall\,v\in\Nodes,\;t\in\Timesteps.
	\label{eq:gen-link}
\end{alignat}
Minimum up-time $\MinUpTime{v} \in \N$ and down-time $\MinDownTime{v} \in \N$ prevent rapid cycling:
\begin{alignat}{3}
	\sum_{t'=t-\MinUpTime{v} + 1}^{t}\zUpt{v}{t'}
	&\le 1 - \zDownt{v}{t}
	&&\quad \forall\,v\in\Nodes,\;t\in\Timesteps,
	\label{eq:gen-minup}
	\\
	\sum_{t'=t-\MinDownTime{v} + 1}^{t}\zDownt{v}{t'}
	&\le 1 - \zUpt{v}{t}
	&&\quad \forall\,v\in\Nodes,\;t\in\Timesteps.
	\label{eq:gen-mindown}
\end{alignat}
\end{subequations}

The numbers of generator startups and shutdowns over the planning horizon can also be limited:
\begin{subequations}\label{eq:gen-action-limits}
	\begin{align}
		\sum_{t\in\Timesteps}\sum_{v\in\Nodes}\zUpt{v}{t}&\leq N^{\mathrm{up}},\\
		\sum_{t\in\Timesteps}\sum_{v\in\Nodes}\zDownt{v}{t}&\leq N^{\mathrm{down}}.
	\end{align}
\end{subequations}

%% ------ Ramps ----------------------------------------
\paragraph{Generation ramps~\textmd{\citep{hedman2010cooptimization}}}
\label{subsec:ramps}
The change in generation of each generator $v \in \Nodes$ between consecutive time steps is limited by ramp rates $\RampDown{v} \in \R$ (downward) and $\RampUp{v} \in \R$ (upward):
\begin{subequations}\label{eq:ramp}
\begin{alignat}{3}
	\Demandt{v}{t-1} - \Demandt{v}{t}
	&\le \RampDown{v}
	&&\quad \forall\,v\in\Nodes,\;t\in\Timesteps,
	\label{eq:ramp-down}\\
	\Demandt{v}{t} - \Demandt{v}{t-1}
	&\le \RampUp{v}
	&&\quad \forall\,v\in\Nodes,\;t\in\Timesteps.
	\label{eq:ramp-up}
\end{alignat}
\end{subequations}

%% ------ Must-run groups ------------------------------
\paragraph{Must-run generator groups}
\label{subsec:must-run}
The set $\MustRunGroups \subseteq 2^\Nodes$ collects groups of generators that must maintain a minimum and maximum aggregate output.
Each group $G^\text{mr}\in\MustRunGroups$ must jointly produce within prescribed bounds:
\begin{alignat}{3}
	\MustRunMint{t}
	\;\le\;
	\sum_{v\in G^\text{mr}}\Demandt{v}{t}
	&\le \MustRunMaxt{t}
	&&\quad \forall\,G^\text{mr}\in\MustRunGroups,\;t\in\Timesteps.
	\label{eq:must-run}
\end{alignat}

%% ------ Transformer groups ---------------------------
\paragraph{Transformer groups}
\label{subsec:trafo-groups}
The set $\TransGroups \subseteq \Arcs\!\times\!\Arcs$ contains pairs of parallel-configured transformers.
They must maintain equal tap positions to avoid circulating currents:
\begin{alignat}{3}
	\Tapt{u}{v}{t}
	&= \Tapt{u'}{v'}{t}
	&&\quad \forall\,\bigl((u,v),(u',v')\bigr)\in\TransGroups,\;t\in\Timesteps.
	\label{eq:trafo-group}
\end{alignat}

%% ------ Angle limits on and off ----------------------
\paragraph{Angle limits in off states \& coupling angle limits~\textmd{\citep{khodaei2012voltage,khodaei2010transmission}}}
\label{subsec:angle-on-off}
\label{subsec:coupling-angle}
In \Cref{eq:angle-diff}, angle differences are constrained only when the arc is on.
Allowing unconstrained off-state angles can produce topologies that are AC-infeasible (see \Cref{subsubsec:feasibility-practice}).
Additionally, coupling can only be performed between buses with low voltage angle difference.
We consolidate these considerations under the bound $\PotDifMaxOff{u}{v} \in \R$, $\PotDifMaxOff{u}{v} \geq \PotDifMax{u}{v}$, on all arcs $(u,v) \in \Arcs\setminus\C$, which is imposed regardless of switching state:
\begin{subequations}\label{eq:angle-off}
\begin{alignat}{3}
	-\PotDifMaxOff{u}{v}
	\;\le\;
	\PotCt{u}{t}-\PotCt{v}{t}-\Shiftt{u}{v}{t}
	&\le
	\PotDifMaxOff{u}{v}
	&&\quad \forall\,(u,v)\in\LArcs\setminus\C\qtime,
	\label{eq:angle-off-branch}\\
	-\PotDifMaxOff{u}{v}
	\;\le\;
	\PotCt{u}{t}-\PotCt{v}{t}
	&\le
	\PotDifMaxOff{u}{v}
	&&\quad \forall\,(u,v)\in\SArcs\setminus\C\qtime.
	\label{eq:angle-off-substation}
\end{alignat}
\end{subequations}
For couplers, if a coupling action is indeed planned in the next time step, a tighter bound~$\PotDifMaxCplAt{(u,v)} \in \R$, $\PotDifMaxCplAt{(u,v)} \leq \PotDifMaxOff{u}{v}$, applies:
\begin{alignat}{3}
	-(1- \zCouplet{u}{v}{t+1}) \PotDifMaxOff{u}{v} - \PotDifMaxCplAt{(u,v)}
	\;\le\;
	\PotCt{u}{t} - \PotCt{v}{t}
	&\le
	(1 - \zCouplet{u}{v}{t+1}) \PotDifMaxOff{u}{v} + \PotDifMaxCplAt{(u,v)}
	&&\quad \forall\,(u,v)\in\Couplers\setminus\C,\;t\in\Timesteps.
	\label{eq:couple-angle-tight}
\end{alignat}

%% ------ Different ratings ----------------------------
\paragraph{Differentiated thermal ratings~\textmd{\citep{hedman2009optimal}}}
\label{subsec:diff-ratings}
For selected arcs, the contingency-case capacity may differ from the base-case value.
This is implemented by replacing $\FlowMaxOver{u}{v}, \FlowMax{u}{v}$ in Equations~\eqref{eq:thermal_over}--\eqref{eq:thermal_sw} with contingency-specific parameters~$\tilde{\FlowParSym}^{\C}_{(u,v)}, \bar{\FlowParSym}^{\C}_{(u,v)} \in \R$ for each such contingency $\C \in \aC$ and arc $(u,v) \in \Arcs$.

%% ------ Branch utilization symmetry ------------------
\paragraph{Branch utilization symmetry}
\label{subsec:branch-sym}
Parallel arcs $\left(a,a'\right)\in\ParArcs$, $\ParArcs \subseteq \Arcs \times \Arcs$, should carry similar load to avoid uneven component aging.
The difference in relative base-case utilization is captured by $\FlowRelDift{a}{a'}{t} \in \R_{\geq 0}$ and cannot exceed $\FlowRelDifMax{a}{a'} \in \R_{\geq 0}$ while both arcs are switched on:
\begin{subequations}\label{eq:sym}
\begin{alignat}{3}
	\frac{\FlowSym_{a}^{t}}{\FlowMaxSym_a} - \frac{\FlowSym_{a'}^{t}}{\FlowMaxSym_{a'}}
	&\le \FlowRelDift{a}{a'}{t} + \bigl(2 - \zSym_{a}^t - \zSym_{a'}^t\bigr)
	&&\quad \forall\,\left(a,a'\right)\in\ParArcs\qtime,
	\label{eq:sym-def1}
	\\
	\frac{\FlowSym_{a'}^{t}}{\FlowMaxSym_{a'}} - \frac{\FlowSym_{a}^{t}}{\FlowMaxSym_a}
	&\le \FlowRelDift{a}{a'}{t} + \bigl(2 - \zSym_{a}^t - \zSym_{a'}^t\bigr)
	&&\quad \forall\,\left(a,a'\right)\in\ParArcs\qtime,
	\label{eq:sym-def2}
	\\
	\FlowRelDift{a}{a'}{t}
	&\le \FlowRelDifMax{a}{a'}
	&&\quad \forall\,\left(a,a'\right)\in\ParArcs\qtime.
	\label{eq:sym-bound}
\end{alignat}
\end{subequations}
Additionally, parallel arcs originating at the same substation must connect to different busbars.
For all parallel arcs $\left((u,v),(u', v')\right)\in\ParArcs$ where, say, $u$ and $u'$ are breaker nodes of the same substation $\Sub \in \Subs$, we add the constraints
\begin{equation}
	\zt{u}{w}{t} + \zt{u'}{w}{t} \leq 1 \quad \forall\, w \in\BusbarsIn{\Sub}\qtime.
	\label{eq:sym-busbar}
\end{equation}

%% ------ Generation symmetry --------------------------
\paragraph{Generation distribution symmetry}
\label{subsec:gen-dist-sym}
Let $\GenNodesIn{\Sub} \subseteq \Nodes$ denote the generator nodes in substation~$\Sub \in \Subs$.
At each substation, generation capacity should be spread evenly across busbars.
With the substation-internal switching variable~$\zt{u}{v}{t}\in\B$ assigning generator node~$v$ to busbar~$u$, the deviation from the ideal equal split is captured by~$\GenDistt{u}{t}\in\R$:
\begin{subequations}\label{eq:gen-dist}
\begin{alignat}{3}
	\sum_{v\in\GenNodesIn{\Sub}}
	\DemandMaxt{\!v}{t}\,\zt{u}{v}{t}
	-
	\frac{1}{|\BusbarsIn{\Sub}|}
	\sum_{v\in\GenNodesIn{\Sub}}\DemandMaxt{\!v}{t}
	&\le \GenDistt{u}{t}
	&&\quad\forall\,\Sub\in\Subs,\;u\in\BusbarsIn{\Sub}\qtime,
	\label{eq:gen-dist1}
	\\
	\frac{1}{|\BusbarsIn{\Sub}|}
	\sum_{v\in\GenNodesIn{\Sub}}\DemandMaxt{\!v}{t}
	-
	\sum_{v\in\GenNodesIn{\Sub}}
	\DemandMaxt{\!v}{t}\,\zt{u}{v}{t}
	&\le \GenDistt{u}{t}
	&&\quad\forall\,\Sub\in\Subs,\;u\in\BusbarsIn{\Sub}\qtime.
	\label{eq:gen-dist2}
\end{alignat}
\end{subequations}

%% ------ Emergency flexibility ------------------------
\paragraph{Emergency flexibility}
\label{subsec:flex}
Generators should be able to restore all branch utilizations below a threshold $\FlexThresh\in[0,1]$.
A \emph{flexibility scenario} is added with its own generation variables $\Demandflext{v}{t} \in \R$, and, as with contingency scenarios, flow $\Flowflext{u}{v}{t} \in \R$, $(u,v) \in \Arcs$ and voltage angle variables $\Potflext{u}{t} \in \R$, $u \in \Nodes$.
These need to satisfy a tightened thermal limit:
\begin{alignat}{3}
	-\zt{u}{v}{t}\,\FlexThresh\,\FlowMax{u}{v}
	\;\le\; \Flowflext{u}{v}{t}
	&\le
	\zt{u}{v}{t}\,\FlexThresh\,\FlowMax{u}{v}
	&&\quad \forall\,(u,v)\in\Arcs\qtime.
	\label{eq:flex}
\end{alignat}
Otherwise, the scenario does not differ from contingency scenarios.

%% ------ Extended objective ---------------------------
\paragraph{Extended objective}
\label{subsec:extended-obj}
With the redispatch, switching-change, coupling, and symmetry variables above, the objective~\eqref{eq:objective} is extended to
\begin{multline}\label{eq:objective-ext}
	\min\;\;
	\sumtime
	\biggl(
	\underbrace{
	\sum_{v\in\Nodes}
	\DemandCostst{v}{t}\,\Demandt{v}{t}
	}_{\text{production cost}}
	\;+\;
	\underbrace{
	\sum_{v\in\Nodes}
	\bigl(\CostRdP{v}\,\DemandPlust{v}{t}
			+\CostRdM{v}\,\DemandMinust{v}{t}\bigr)
	}_{\text{redispatch cost}}
	+\;
	\underbrace{
	\sum_{\C \in \aC}
	\sum_{(u,v)\in\Arcs}
	\OverCost\,\FlowOverCt{u}{v}{t}
	}_{\text{overload penalty}}
	\\
	+\;
	\underbrace{
	\sum_{(u,v)\in\Arcs}
	\CostCh{u}{v}\,\zChanget{u}{v}{t}
	}_{\text{switching cost}}
	\;+\;
	\underbrace{
	\sum_{(u,v)\in\Couplers}
	\bigl(\CostCp{u}{v}\,\zCouplet{u}{v}{t}
			+\CostSp{u}{v}\,\zSplitt{u}{v}{t}\bigr)
	}_{\text{coupling/splitting cost}}
	\;+\;
	\underbrace{
	\sum_{\Sub \in \Subs}
	\sum_{u \in \BusbarsIn{\Sub}}
	\GenDistt{u}{t}
	}_{\text{symmetry penalty}}
	\biggr).
\end{multline}

\paragraph{Implementation scope}
The formulations above include symmetric branch utilization, generation-distribution symmetry, and emergency flexibility, which are not implemented in the solver described here.
AC-only constraints are outside the DC optimization model; selected solutions are checked separately using AC power flow calculations (\Cref{sec:practice}).

\section{Methodology}\label{sec:methodology}
This section describes the algorithmic components to exactly solve the previously introduced MILP models.
Our focus lies on the core problem, \Cref{subsec:dc-ots}, \ie, a DC-SCOTS with substation reconfiguration, transformers, and generator (``unit'') commitment.
We also discuss handling of the problem extensions (\Cref{subsec:extensions}) to which the methodology does not naturally extend.
The superscript for the fixed time step $t \in \Timesteps$ is omitted from the notation in this section.

Our approach is a branch-and-Benders-cut algorithm that decomposes the problem naturally along the contingency structure, with each contingency represented by a subproblem.
As this procedure allows for individual consideration of contingencies, it avoids the potentially intractable memory requirements caused by replicating the power-flow constraints for every contingency and time step in the compact model.
\Cref{subsec:bbc} provides additional details.
The rest of the section roughly follows a bottom-up (SP-MP) structure.
\Cref{subsec:les} introduces a specialized subproblem algorithm.
Ranking and selecting subproblems reduces the number of SPs that are solved, while simultaneously limiting bloating of the MP, and is described in \Cref{subsec:subproblem-selection}.
The latter is further improved by cut filtering, \Cref{subsec:cut-filtering}.
Partial decomposition, in \Cref{subsec:partial-decomposition}, moves certain SPs into the MP, reducing the number of iterations until Benders convergence.
We outline a connectivity cut separation procedure in \Cref{subsec:connectivity-method}.
Finally, \Cref{subsec:sequential-switching-heuristic} describes the sequential switching heuristic and its use as a warmstart or parallel heuristic for the exact solve to improve primal solution quality.

\subsection{Branch-and-Benders-Cut}
\label{subsec:bbc}

The original problem is decomposed into a master problem (MP) and one subproblem (SP) per contingency $\C\in\aC$  (and time step).
The MP contains the topology, transformer, commitment, and redispatch variables and their associated Constraints~\eqref{eq:tap}, \eqref{eq:breaker}--\eqref{eq:gen-bounds}, \eqref{eq:rd}--\eqref{eq:default}, \eqref{eq:gen-updown}--\eqref{eq:trafo-group}, \eqref{eq:sym-busbar}, and~\eqref{eq:gen-dist}.
The remaining power-flow and phase-angle constraints are considered in the subproblems.
One auxiliary variable $\eta_\C\geq 0$ is added to the MP for each contingency $\C\in\aC$.
The objective is altered to
\begin{equation}
	\min\;\text{(Equation \eqref{eq:objective-ext} without overload or symmetry penalties)}
		\;+\;\sum_{\C\in\aC}\eta_\C.
	\label{eq:mp-obj}
\end{equation}
The MP is solved via branch-and-bound.
At each node with an integer-feasible (or fractional) solution, the master decisions $\zFixedSym \in \B^\Arcs$, $\DemandFixedSym \in \R^\Nodes$, $\ShiftFixedSym \in \R^\LArcs$, $\etaFixedSym \in \R^\aC$ are first checked for connectivity, then passed to the subproblems for separation.
We maintain one MP and reuse one SP template per worker across contingencies by updating the contingency-dependent constraints.
Whenever the master solution changes, we also update the fixed master decisions in each template.

\paragraph{Subproblem}
Given a fixed master solution $\zFixedSym$, $\DemandFixedSym$, $\ShiftFixedSym$, $\etaFixedSym$, we denote the subproblem corresponding to contingency $\C \in \aC$ by $\subP$.
It checks whether the corresponding contingency power flow is feasible and if the auxiliary variable $\eta_\C$ does not underestimate the subproblem's cost, and returns a cut otherwise.
Contingency branches in $\C$ are removed and the master variable values treated as parameters.
Denoting $\Pot{(u,v)}^{\C,\ShiftSym} \coloneq \PotC{u}\!-\!\PotC{v}\!-\!\ShiftFixed{u}{v}$, the subproblem reads:
% fixed-width box for the dual annotation so the \forall column aligns
\newcommand{\spdual}[1]{\;\mathmakebox[\widthof{$[\,\DualOverCLbSym,\DualOverCUbSym\,]$}][l]{#1}}
\begin{mini!}|s|
	{}
	{\sum_{(u,v)\in\Arcs\setminus\C}\OverCost\,\FlowOverC{u}{v}{\C}\label{eq:sp-obj}}
	{\label{eq:sp}}{}
	%% power balance
	\addConstraint{\sum_{(v,u)\in\delta^-(u) \setminus\C}\!\FlowC{v}{u} - \sum_{(u,v)\in\delta^+(u) \setminus\C}\!\FlowC{u}{v}}
		{= \StatD{u} - \DemandFixed{u}}
		{\spdual{[\DualFCCSym]}\quad \forall\, u\in\Nodes\label{eq:sp-pb}}
	%% Ohm's law
	\addConstraint{\left(\Pot{(u,v)}^{\C,\ShiftSym} - \bigl(1\!-\!\zFixed{u}{v}\bigr) \PotDifMaxOff{u}{v} \right) \Sus{u}{v}}
		{\le \FlowC{u}{v}}
		{\spdual{[\DualPECLbSym]}\quad \forall\,(u,v)\!\in\!\LArcs\!\setminus\!\C\label{eq:sp-ohm-lb}}
	\addConstraint{\left(\Pot{(u,v)}^{\C,\ShiftSym} + \bigl(1\!-\!\zFixed{u}{v}\bigr) \PotDifMaxOff{u}{v} \right) \Sus{u}{v}}
		{\ge \FlowC{u}{v}}
		{\spdual{[\DualPECUbSym]}\quad \forall\,(u,v)\!\in\!\LArcs\!\setminus\!\C\label{eq:sp-ohm-ub}}
	%% thermal limits with overload
	\addConstraint{-\bigl(\FlowMaxOver{u}{v}+\FlowOverC{u}{v}{\C}\bigr) \;\le\; \FlowC{u}{v}}
		{\le \FlowMaxOver{u}{v}+\FlowOverC{u}{v}{\C}}
		{\spdual{[\DualOverCLbSym,\DualOverCUbSym]}\quad \forall\,(u,v)\in\LArcs\!\setminus\!\C\label{eq:sp-thermal_over}}
	%% switch coupling
	\addConstraint{-\zFixed{u}{v} \FlowMax{u}{v} \;\le\;\FlowC{u}{v}}
		{\le \zFixed{u}{v} \FlowMax{u}{v}}
		{\spdual{[\DualSOffCLbSym,\DualSOffCUbSym]}\quad \forall\,(u,v)\!\in\!\Arcs\!\setminus\!\C\label{eq:sp-switch}}
	%% voltage angle difference
	\addConstraint{-(1\!-\!\zFixed{u}{v})\PotDifMaxOff{u}{v} \!-\! \PotDifMax{u}{v}}
		{\le \Pot{(u,v)}^{\C,\ShiftSym}}
		{\spdual{[\DualSOnCLbSym]}\quad \forall\,(u,v)\!\in\!\LArcs\!\setminus\!\C\label{eq:sp-vad-lb}}
	\addConstraint{(1\!-\!\zFixed{u}{v})\PotDifMaxOff{u}{v} \!+ \PotDifMax{u}{v}}
		{\ge \Pot{(u,v)}^{\C,\ShiftSym}}
		{\spdual{[\DualSOnCUbSym]}\quad \forall\,(u,v)\!\in\!\LArcs\!\setminus\!\C\label{eq:sp-vad-ub}}
	%% substation-internal branch coupling
	\addConstraint{-(1 - \zFixed{u}{v}) \PotDifMaxOff{u}{v}}
		{\leq \PotC{u} - \PotC{v}}
		{\spdual{[\DualSOnCLbSym]}\quad \forall\,(u,v)\!\in\!\SArcs\label{eq:sp-aux-lb}}
	\addConstraint{(1 - \zFixed{u}{v}) \PotDifMaxOff{u}{v}}
		{\geq \PotC{u} - \PotC{v}}
		{\spdual{[\DualSOnCUbSym]}\quad \forall\,(u,v)\!\in\!\SArcs\label{eq:sp-aux-ub}}
	%% angle off limits
	\addConstraint{-\PotDifMaxOff{u}{v} \;\le\; \Pot{(u,v)}^{\C,\ShiftSym}}
		{\le \PotDifMaxOff{u}{v}}
		{\spdual{[\DualVadOffCLbSym,\DualVadOffCUbSym]}\quad \forall\,(u,v)\in\LArcs\!\setminus\!\C\label{eq:sp-angle-off-branch}}
	\addConstraint{-\PotDifMaxOff{u}{v} \;\le\; \PotC{u}-\PotC{v}}
		{\le \PotDifMaxOff{u}{v}}
		{\spdual{[\DualVadOffCLbSym,\DualVadOffCUbSym]}\quad \forall\,(u,v)\in\SArcs.\label{eq:sp-angle-off-substation}}
\end{mini!}
We note the dual variables, or ``Farkas multipliers,'' corresponding to the respective constraints in square brackets.

\paragraph{Cut generation}
Three cases arise after solving $\subP$:
\begin{enumerate}
\item
	\textbf{Optimality Cut.}
	If the subproblem is feasible and its optimal value $\subP^{*} \in \R$ exceeds the current estimate, i.e., $\subP^{*} > \etaFixedSym_\C$, a corresponding dual solution yields an optimality cut. Denoting the dual solution values by an asterisk ($*$), the cut reads
	\begin{subequations}\label{eq:opt-cut-h}
	\begin{equation}
		\eta_\C \ge
		\sum_{u\in\Nodes}
			(\StatD{u}-\Demand{u})\,\DualFCCOpt{u}
		\;+\!\!\!
		\sum_{(u,v)\in\Arcs\setminus\C}
			h^{\C*}_{(u,v)},
		\label{eq:opt-cut}
	\end{equation}
	where each term $h^{\C*}_{(u,v)}$ collects the dual contributions that depend on arc $(u,v) \in \Arcs \setminus \C$ via the big-$M$ right-hand sides of \eqref{eq:sp-ohm-lb}--\eqref{eq:sp-angle-off-substation}. For any branch $(u,v) \in \LArcs \setminus \C$,
	\begin{equation}
	\begin{split}
		h^{\C*}_{(u,v)} ={}& \left((1\!-\!\z{u}{v})\Sus{u}{v}\PotDifMaxOff{u}{v} + \Sus{u}{v}\Shift{u}{v}\right)\,\DualPECLbOpt{u}{v} \\
		&+\left((1\!-\!\z{u}{v})\Sus{u}{v}\PotDifMaxOff{u}{v} - \Sus{u}{v}\Shift{u}{v}\right)\,\DualPECUbOpt{u}{v} \\
		&+\left(\FlowMaxOver{u}{v}\right)\left(\DualOverCLbOpt{u}{v}+\DualOverCUbOpt{u}{v}\right)
		 +\left(\z{u}{v} \FlowMax{u}{v}\right) \left(\DualSOffCLbOpt{u}{v}+\DualSOffCUbOpt{u}{v}\right) \\
		&+\left((1\!-\!\z{u}{v})\PotDifMaxOff{u}{v} + \PotDifMax{u}{v} - \Shift{u}{v}\right)\DualSOnCLbOpt{u}{v}\\
		&+\left((1\!-\!\z{u}{v})\PotDifMaxOff{u}{v} + \PotDifMax{u}{v} + \Shift{u}{v}\right)\DualSOnCUbOpt{u}{v}\\
		&+\left(\PotDifMaxOff{u}{v} - \Shift{u}{v}\right)\DualVadOffCLbOpt{u}{v}
		 +\left(\PotDifMaxOff{u}{v} + \Shift{u}{v}\right)\DualVadOffCUbOpt{u}{v},
	\end{split}
		\label{eq:opt-cut-h-branch}
	\end{equation}
	and for any substation arc $(u,v) \in \SArcs$,
	\begin{equation}
	\begin{split}
		h^{\C*}_{(u,v)} ={}&\left(\z{u}{v} \FlowMax{u}{v}\right) \left(\DualSOffCLbOpt{u}{v}+\DualSOffCUbOpt{u}{v}\right) 
		 +\left((1\!-\!\z{u}{v})\PotDifMaxOff{u}{v}\right)\left(\DualSOnCLbOpt{u}{v}+\DualSOnCUbOpt{u}{v}\right)\\
		&+\left(\PotDifMaxOff{u}{v}\right) \left(\DualVadOffCLbOpt{u}{v}+\DualVadOffCUbOpt{u}{v}\right).
	\end{split}
		\label{eq:opt-cut-h-substation}
	\end{equation}
	\end{subequations}

\item
	\textbf{Feasibility cut.}
	If the subproblem is infeasible
	a dual ray with positive objective value replaces the dual solution in~\eqref{eq:opt-cut}, and the left-hand side becomes~$0$:
	\begin{equation}
	0 \;\ge\;
	\sum_{u\in\Nodes}
		(\StatD{u}-\Demand{u})\,\DualFCCOpt{u}
	+
	\sum_{(u,v)\in\Arcs\setminus\C}
		h^{\C*}_{(u,v)}.
	\label{eq:feas-cut}
	\end{equation}

\item
	\textbf{No cut.}
	If $\subP^{*}\le\etaFixedSym_\C$, the contingency is already covered by the current master solution and no cut is needed.
\end{enumerate}

\subsection{Solving Subproblems as Linear Equation Systems}
\label{subsec:les}

The approach summarized in this subsection is developed in detail in~\textcite{donkiewicz2026potflow} for the broader class of affine potential-based flow problems; we discuss its application to our subproblems.
The approach consists of two phases: first, we determine a primal solution to check for constraint violations, then we determine a dual solution or direction of unboundedness to construct cuts.

\paragraph{Primal solution}
To obtain a primal solution, we use the fact that the system of linear equations (SLE) corresponding to a power flow on a fixed network has a unique solution \citep{birkhoff1956potential}.
The SLE is derived from Equations \eqref{eq:sp-pb}--\eqref{eq:sp-angle-off-substation} as follows.
Given a fixed topology $\zFixed{u}{v} \in \B^\Arcs$, the set of \emph{active arcs} $\ArcsOutz \coloneq \{(u,v)\in\Arcs \mid \zFixed{u}{v}=1\}\setminus\C$ is fully known.
We correspondingly denote the sets of branches $\LArcsCC \coloneq \LArcs \cap \ArcsOutz$ and substation arcs $\SArcsz \coloneq \SArcs \cap \ArcsOutz$.
The flow on the arcs outside $\ArcsOutz$ is zero, and they do not further influence the power flow on other arcs, and can thus be disregarded.
The big-$M$ terms in the remaining constraints can then be resolved for active arcs, \ie, we can remove the terms for Constraints \eqref{eq:sp-ohm-lb}, \eqref{eq:sp-ohm-ub}, \eqref{eq:sp-vad-lb} -- \eqref{eq:sp-aux-ub}, and remove the Constraints \eqref{eq:sp-switch} for substation arcs.
This makes the substation arcs redundant, they only force the voltage angles of incident nodes to be equal.
We ``contract'' them by replacing each connected component from the graph containing only active substation arcs $(\Nodes, \SArcsz)$ by a single node which is connected to all active branches incident to nodes in the component.
The contracted graph contains no more substation arcs.

Without the big-$M$ terms, the Ohm's law constraints~\eqref{eq:sp-ohm-ub}--\eqref{eq:sp-ohm-lb} turn into the equality $\FlowC{u}{v}=\Sus{u}{v}(\PotC{u}-\PotC{v}-\ShiftFixed{u}{v})$ for all $(u,v)\in\LArcsCC$.
We substitute the flow variables in the power-balance constraint~\eqref{eq:sp-pb}, resulting in the system of linear equations $\lesP$:
\begin{multline}
	  \sum_{(v,u)\in\delta^-(u) \setminus\C}\!(\PotC{v}-\PotC{u})\Sus{v}{u}
	- \sum_{(u,v)\in\delta^+(u) \setminus\C}\!(\PotC{u}-\PotC{v})\Sus{u}{v}\\
	= \StatD{u} - \DemandFixed{u} + \sum_{(v,u)\in\delta^-(u) \setminus\C}\!\ShiftFixed{v}{u}\Sus{v}{u} - \sum_{(u,v)\in\delta^+(u) \setminus\C}\!\ShiftFixed{u}{v}\Sus{u}{v} \quad \forall u \in \Nodes
	\label{eq:les-primal}
\end{multline} 
The coefficient matrix is the weighted Laplacian of the contracted graph.
The system $\lesP$ yields a solution that is unique in voltage angle differences and flows \citep{birkhoff1956potential}.
We can thus use its solutions to check for violations of the remaining constraints \eqref{eq:sp-switch}--\eqref{eq:sp-vad-ub}, \eqref{eq:sp-angle-off-branch}--\eqref{eq:sp-angle-off-substation}.
If there are violations, the SP is infeasible and a feasibility cut needs to be generated.
Otherwise, the optimal overload of branch $(u,v) \in \LArcsCC$ is given by $\FlowOverCOpt{u}{v}=\max\{|\FlowCOpt{u}{v}|-\FlowMaxOver{u}{v},\,0\}$, and the objective value is $\subP^{*}=\sum_{(u,v)}\OverCost\,\FlowOverCOpt{u}{v}$.
If it is larger than $\etaFixedSym_\C$, an optimality cut needs to be generated.

\paragraph{Optimality and feasibility cuts}
To construct an optimality cut, we require a feasible dual solution.
By using a basis obtained from the primal solution computed in the first step, we construct another SLE, this time in the $\DualFCCSym$ dual variables.
Again, buses and arcs are those of the contracted graph. 
For a given ordering of buses, the coefficient matrix is equal to that of the primal system, only the right-hand side changes: 
\begin{multline}
	  \sum_{(v,u)\in\delta^-(u) \setminus\C}\!(\DualFCC{v}-\DualFCC{u})\Sus{v}{u}
	- \sum_{(u,v)\in\delta^+(u) \setminus\C}\!(\DualFCC{u}-\DualFCC{v})\Sus{u}{v}\\
	= \sum_{(v,u)\in \delta^-(u) \setminus \C} \DualOverCOpt{v}{u} \Sus{v}{u} - \sum_{(u,v)\in \delta^+(u) \setminus \C} \DualOverCOpt{u}{v} \Sus{u}{v} \quad \forall u \in \Nodes,
	\label{eq:les-dual-opt}
\end{multline}
where for each active branch $(u,v) \in \LArcs$, the difference $\DualOverCOpt{u}{v} = \DualOverCLbOpt{u}{v}-\DualOverCUbOpt{u}{v}$ of the dual variables of the overload constraints \eqref{eq:sp-thermal_over}, $\DualOverCLb{u}{v}, \DualOverCUb{u}{v}$, is fixed based on whether and in which direction the branch is overloaded:
\begin{align*}
	\DualOverCOpt{u}{v}   &= \begin{cases}
                            0, & -\FlowMaxOver{u}{v} \leq \FlowCOpt{u}{v} \leq \FlowMaxOver{u}{v} \\
                            -\OverCost, & \FlowCOpt{u}{v} < -\FlowMaxOver{u}{v} \\
                            \OverCost, & \FlowCOpt{u}{v} > \FlowMaxOver{u}{v}.
                          \end{cases}
\end{align*}
The remaining dual solution values are then either zero or can be determined directly from $\DualFCCOptSym$.
All dual solutions determined this way are optimal.

To construct a feasibility cut, we similarly require a dual direction with positive objective value.
The coefficient matrix of the SLE here is an expanded version of that of optimality cuts, with an additional column to represent an arc where a violation occurs, right-hand side $0$, and an additional row corresponding to the dual objective function, right-hand side $1$.

\paragraph{Factorization updates}
The coefficient matrices between Benders iterations, contingencies, and the primal and dual SLEs of the same subproblem are very similar or equal.
We can make use of this by updating the $LU$ factorization of the coefficient matrix from subproblem to subproblem.
After computing an initial factorization, at most one rank-$1$ update is required per arc that differs between subproblems.
This includes both the topology given by the master solution and the contingency corresponding to the subproblem.
Roughly speaking, solving and computing a factorization for an SLE is cubic in the number of buses, whereas updating an existing, sufficiently similar factorization, and then solving it is quadratic.
This reduces the overall complexity of solving all subproblems of one Benders iteration from $\mathcal{O}(|\Arcs| \cdot |\Nodes|^3)$ to $\mathcal{O}(|\Arcs| \cdot |\Nodes|^2)$.
We refer to this approach as SLE-U.

\subsection{Subproblem Selection}
\label{subsec:subproblem-selection}
Due to the large number of oftentimes similar subproblems, the generated cuts can lead the B\&BC algorithm in similar directions.
Additionally, some subproblems do not produce any cuts.
Solving them is only necessary when the master solution is accepted.
Moreover, adding cuts to the master problem can increase the running time needed to solve it.
Solving certain subproblems can thus yield only a small improvement with regard to the number of Benders iterations or branch-and-bound nodes, while incurring unnecessary computational effort for solving the SP and MP.

We hence solve only a subset of all subproblems.
The goal is to obtain a set of cuts, produced by the subproblems, that is overall as beneficial as possible to the remaining B\&BC execution, while minimizing the number of subproblems that are solved.
Two components facilitate this.
Subproblems are sorted at the beginning of each Benders iteration based on their \emph{score}.
The score is based on subproblem metrics collected during the running time.
Subproblems are then solved in order of their score until one of multiple \emph{stopping criteria} is met.
These are generally statically defined at the beginning of the B\&BC execution.

\paragraph{Subproblem scoring}
For each contingency $\C\in\aC$ (and corresponding subproblem), we maintain a score $\spS{\C}{r}\in[0,1]$ composed of the weighted sum of six components.
The first set of these are based on metrics of individual subproblems, with the goal of maximizing the chance that a cut is produced, and the impact of that cut.
Let $r\in\N$ be the number of the current Benders iteration.
Let $\spItSolved{r}{\C} \subseteq \{1,\ldots,r-1\}$ be iterations where (the subproblem corresponding to) $\C$ was solved, and $\spItCut{r}{\C} \subseteq \spItSolved{r}{\C}$ be iterations where solving subproblem $\C$ produced a cut.
For iterations $r' \in \spItCut{r}{\C}$ where an optimality cut was produced, let $\eta_{\C}^*(r') \in \R$ be the subproblem's objective function value, and $\etaFixedSym_{\C}(r') \in \R$ the master solution value of the Benders variable for the subproblem in that iteration.
The \emph{violation} of the subproblem objective function value by the master solution is given by $\Viol{\C}{r'} \coloneq \eta_{\C}^*(r') - \etaFixedSym_{\C}(r')$.
The scores are then as follows:
\begin{itemize}
	\item \emph{Reliability}
	$\spSrel{\C}{r} \coloneq |\spItCut{r}{\C}|/|\spItSolved{r}{\C}|$:
	fraction of solves that produced a cut.
	\item \emph{Total share}
	$\spSshr{\C}{r} \coloneq |\spItCut{r}{\C}|\big/\sum_{\C'\in\aC}|\spItCut{r}{\C'}|$:
	share of all generated cuts.
	\item \emph{Average violation}
	$\spSvio{\C}{r} \coloneq |\spItCut{r}{\C}|^{-1}\sum_{r'\in \spItCut{r}{\C}}\Viol{\C}{r'}$:
	mean overload-cost violation when an optimality cut was generated (subproblem objective function value). Feasibility cuts do not contribute to average violation.
	\item \emph{Staleness}
	$\spSstl{\C}{r} \coloneq (r-1) - \max\{r': r'\in \spItSolved{(r-1)}{\C}\}$:
	number of iterations since $\C$ was last examined.
	\item \emph{Diversity}:
	we calculate a cosine distance-based $k$-medoid clustering~\cite{kaufman1987clustering} of the subproblems, with $k = \spsimk \in \N$, on line-outage distribution factors \citep[cf.][]{guler2007generalized} at the start of the B\&BC algorithm.
	Subproblems are chosen randomly in circular order from the clusters, starting with the medoids.
	The $i$-th subproblem $\C$ chosen this way, $i \in \{1,\ldots,|\aC|\}$, $\C \in \aC$, is assigned diversity score $\spSdiv{\C}{r} \coloneq |\aC| - i$.
	\item \emph{Random component}
	$\spSrnd{\C}{r} \sim \text{Uniform}(0,1)$.
\end{itemize}
Since subproblem relevance might change throughout the solving process, we additionally employ an exponential decay where applicable, effectively multiplying the past score with a decay factor $\spdecay \in [0,1]$ at every iteration.

The final score $\spS{\C}{r}$ is then determined using a weighted combination with the weights $\spWrel, \spWshr, \spWvio, \spWstl, \spWrnd, \spWdiv \ge 0$ balancing the impact of each score.
The subproblems are always sorted in descending order according to $\spS{\C}{r}$. 
Then, they are solved in that order until a stopping criterion is met.

\paragraph{Stopping criteria}
If one or more of the stopping criteria are met, the subproblem solving iteration is stopped, and the cuts found so far are either to be added directly to the master problem, or filtered further (see \Cref{subsec:cut-filtering}).
We call rounds in which all subproblems are solved and all generated cuts are added to the master problem \emph{stabilization rounds}.
They occur every $\spstab \in \N$ rounds.

Let $\spCuts{r} \subseteq \aC$ be the subproblems that were already solved and yielded a cut in this iteration $r$.
Let $\ObjMP{r} \in \R$ be the current master problem objective function value, and $\ObjInc{r} \in \R$ that of the current incumbent.
Outside stabilization rounds ($r \mod \spstab \ne 0$), the iteration may be stopped when $|\spCuts{r}| \ge 1$ and any of the following conditions holds:

\begin{enumerate}
	\item \emph{Absolute cut limit} $|\spCuts{r}| \ge \spmaxn$:
	the number of found cuts exceeds an absolute threshold $\spmaxn \in \N_{\geq 1}$.
	\item \emph{Relative cut limit} $|\spCuts{r}| \ge \spmaxrel |\aC|$:
	the number of found cuts exceeds a fraction $\spmaxrel \in \R_{\geq 0}$ of total subproblems available.
	\item \emph{Score threshold} $\max(\{\spS{\C}{r}: \C \in \aC \setminus \spCuts{r}\}) < \spminscore$:
	the score of the next subproblem to solve is below $\spminscore \in [0,1]$.
	\item \emph{Primal bound threshold} $\left(\ObjMP{r} + \sum_{\C \in \spCuts{r}} \Viol{\C}{r}\right)(1-\sppb) \ge \ObjInc{r}$:
	the current master problem solution value (including cost estimates), corrected by the accrued violations of the already produced cuts, reaches the objective function value of the incumbent up to a relative buffer $\sppb \in [0,1)$. For a feasibility cut, $\Viol{\C}{r}$ denotes the corresponding constraint violation.
\end{enumerate}

Let $\spCutsFeas{r} \subseteq \spCuts{r}$ and $\spCutsOpt{r} \subseteq \spCuts{r}$ be the sets of subproblems that produced a feasibility or optimality cut, respectively.
For improved numerical and convergence stability, at least one of (1)--(4) and one of the following must be satisfied:
\begin{enumerate}\setcounter{enumi}{4}
	\item \emph{Minimum feasibility cuts} $|\spCutsFeas{r}| \ge \spminfeas$:
	the number of feasibility cuts $|\spCutsFeas{r}|$ added is at least $\spminfeas \in \N$.
	\item \emph{Minimum optimality cuts} $|\spCutsOpt{r}| \ge \spminopt$:
	the number of optimality cuts $|\spCutsOpt{r}|$ added is at least $\spminopt \in \N$.
\end{enumerate}

The methodology discussed in this section is the foundation of the non-ML-based methods described in~\textcite{donkiewicz2026subproblem}.

\subsection{Cut Filtering}
\label{subsec:cut-filtering}

In multi-subproblem Benders implementations, the number of cuts added to the master problem can pose a major computational burden.
It can cause the share of time spent solving the master problem to be relatively large, even if the master problem itself is rather simple.

To address this issue, we apply a cut filtering strategy, reducing the number of cuts added to the master problem.
We combine cut filtering with subproblem selection (\Cref{subsec:subproblem-selection}), selecting among the cuts produced by the chosen subproblems.
Instead of using subproblem properties and results from previous iterations to select subproblems and the resulting, yet unknown, cuts, we use properties of the current iteration to filter the already produced cuts.
Within the scope of this paper, we rank cuts by their violation (\ie, the difference of the subproblem objective function and the current master subproblem cost estimate).
Additional filtering strategies include diversity-based filtering and a hybrid of the two.
Outside stabilization rounds, we select the first $\cfn{r}$ cuts by applying the following filtering rules:

\begin{enumerate}
	\item \emph{Absolute number of cuts} $\cfn{r} \le \cfmaxn$:
	only the $\cfmaxn \in \N$ most violated cuts are added to the master problem.
	\item \emph{Relative number of cuts} $\cfn{r} \le \lceil\cfmaxrel |\aC|\rceil$:
	only the $\cfmaxrel$-fraction, $\cfmaxrel \in \R_{\geq 0}$, of most violated cuts are added to the master problem.
	\item \emph{Primal bound threshold} $\sum_{i=1}^{\cfn{r}} \cfSvio{i}{r} \ge \cfpb\,(\ObjInc{r} - \ObjMP{r})$:
	cuts are added in decreasing order of violation until the cumulative violation of the added cuts, $\sum_{i=1}^{\cfn{r}} \cfSvio{i}{r}$, reaches a fraction $\cfpb \in \R_{\geq 0}$ of the current gap $\ObjInc{r} - \ObjMP{r}$ between the incumbent and the master problem solution value.
\end{enumerate}
Further enhancements include adding the discarded cuts as one (weighted) aggregated cut.
Whenever cuts are generated, at least one is retained after filtering.

We refer to \citet{donkiewicz2026cutfilter} for further discussion of cut filtering strategies for related problems.

\subsection{Partial Decomposition}
\label{subsec:partial-decomposition}

Where subproblem selection and cut filtering are intended to reduce the solving time of individual Benders iterations, partial decomposition is intended to reduce the number of iterations.
For this, we retain a certain number of subproblems in the MP.
This can reduce the number of Benders iterations, while increasing the time needed per MP solve.
A potential violation in one of the retained subproblems disincentivizes solutions that might incur violations in similar, remaining subproblems as well.
As our subproblems only differ in the low number of contingency branches, few retained subproblems can already significantly reduce the number of subproblems with violations.

We use branch relevance scores to determine which branch contingencies (and their associated subproblems) to include in the master problem.
In the following, let $\pdflow{(u,v)} \in \R$ denote the power flow on branch $(u,v)$ in a reference solution (in our case, an optimal OPF solution).
For a congestion threshold $\lambda\in(0,1]$, the set of \emph{congested branches} is
\begin{equation*}
  \pdcong
  = \bigl\{(u,v)\in\Arcs : |\pdflow{(u,v)}| \ge \lambda\cdot\FlowMax{u}{v}\bigr\}.
\end{equation*}

We employ the following scores:

\begin{itemize}
	\item Congestion proximity (\textsc{CongProx}):
	Let $\HopDist{a}{a'} \in \N$ be the hop distance between branches $a \in \Arcs$ and $a' \in \Arcs$, \ie, the length of a shortest path connecting them.
	For a hop radius $\pdcphop\in\N$, the score is
	\begin{equation}
		\pdScp{a}
		\coloneq \max\Bigl(0,\;\max_{\substack{a'\in\pdcong\\\HopDist{a}{a'}\le \pdcphop}} \tfrac{1}{1+\HopDist{a}{a'}}\Bigr),
		\label{eq:score-congprox}
	\end{equation}
	for any branch $a \in \Arcs$, where the inner maximum equals $0$ if the set is empty.
	Branches not within $\pdcphop$ hops of any congested branch receive score~$0$.
	Any congested branch $a\in\pdcong$ has $\HopDist{a}{a}=0$ and thus receives the maximum score~$1$.

	\item Branch switching sensitivity (\textsc{Sens}, \citet{liu2012heuristic}):
	The score measures the maximum relative flow reduction on congested branches when switching off branch $a \in \Arcs$, based on \emph{line outage distribution factors} (LODFs).
	The LODF $\LODF{a'}{a} \in \R$ of branches $a', a \in \Arcs$ is the change in power flow on $a'$, relative to the pre-outage flow on $a$, after $a$ is removed from the network.
	In a DC power flow model, the LODF is computed from the power-transfer distribution factor (PTDF) which in turn can be calculated using the reduced node-susceptance matrix \citep[cf.][]{guler2007generalized}.
	The sensitivity score is given by
	\begin{equation}
		\pdSsns{a}
		= -\min_{\substack{a'\in\pdcong\\\pdflow{a'} \ne 0}}
		  \frac{|\pdflow{a'} + \LODF{a'}{a}\,\pdflow{a}| - |\pdflow{a'}|}{|\pdflow{a'}|}.
		\label{eq:score-sens}
	\end{equation}
	A negative minimum (\ie, a flow reduction on all congested branches) yields a positive score.
	The less ``already overloaded'' branches are additionally overloaded, the higher the score.

	\item Subproblem cost (\textsc{Violation}):
	The score equals the total estimated overload cost by contingencies that contain branch~$a \in \Arcs$:
	\begin{equation}
		\pdSvio{a} = \sum_{\C\in\aC:\,a\in\C} \etaFixedSym_\C,
		\label{eq:score-eta}
	\end{equation}
	where $\etaFixedSym_\C \ge 0$ is the current master problem cost estimate for contingency $\C$.
\end{itemize}

The final score $\pdS{\C} \in [0, 1]$ of contingency $\C \in \aC$ is a weighted combination of the normalized individual scores, with weights $\pdWcp, \pdWsns,  \pdWvio \in [0, 1]$ for congestion proximity, sensitivity, and violation.
Subproblems are then ranked according to those scores, with higher ranked scores being more likely to be included in the MP.
Optionally, the base case, \ie, the empty contingency $\emptyset$, can be prioritized to be included first.

We then partition the set of contingencies $\aC$ into $\pdinc$ and $\pddec$, where the former is included in the master problem and the latter is decomposed as before.
The number of included contingencies is limited by an absolute number of contingencies, $|\pdinc| \le \pdmaxn$, $\pdmaxn \in \N$, or a fraction of the number of contingencies, $|\pdinc| \le \pdmaxrel \cdot |\aC|$, $\pdmaxrel \in [0,1]$.
Then,
\begin{equation}
	\pdinc = \argmax_{S \subseteq \aC} \{\sum_{\C \in S} \pdS{\C} \mid |S| \leq \min\{\pdmaxrel \cdot |\aC|, \pdmaxn\} \},\qquad \pddec = \aC\setminus\pdinc.
\end{equation}

As an alternative to the score, we also implement a diversity-based strategy (\textsc{Div}).
It uses the same distance measure and is equally motivated as for subproblem selection, \ie, the distance between two contingencies is equal to the cosine distance of their LODF vectors, as contingencies with similar LODF vectors might exhibit similar behavior.
The contingencies to be included, optionally next to the base case, are then the medoids resulting from a $k$-medoid clustering, with $k = \min\{\pdmaxrel \cdot |\aC|, \pdmaxn\}$ equal to the desired number of subproblems.

\subsection{Network Connectivity}
\label{subsec:connectivity-method}
As discussed for Constraints \eqref{eq:connected}, the network has to be islanding-free.
In a preprocessing step, we fix all branches on that would already individually disconnect in conjunction with a contingency.
In addition, since the number of possible island configurations is exponential in the number of arcs, and the constraints themselves thus impractical, connectivity is enforced via a separation algorithm, placed before the Benders iteration.
To detect a violation, we count the number of islands, \ie, components in the graph induced by the switched-on arcs that have at least one node with non-zero demand or generation not fixed to zero.
If there is already more than one island, there trivially is also more than one island if the branches of any contingency $\C \in \aC$ are removed as well.
Otherwise, we then search for a contingency that would split the single island into two islands by computing its bridges.
For the first such contingency $\C \in \aC$ we find, we add the connectivity cut
\begin{equation}
	\sum_{(u,v)\in\Arcs\setminus(\ArcsOutz\cup\C)} \z{u}{v} \;\ge\; 1,
	\label{eq:conn-cut}
\end{equation}
which forces at least one of the previously switched-off arcs to be switched on.

\subsection{Sequential Switching Heuristic}
\label{subsec:sequential-switching-heuristic}

Under operational time limits, we are particularly interested in good primal solutions.
For any feasible instance, at least one feasible solution should be found.
We initially restrict the set of available topologies in order to obtain feasible, but potentially suboptimal solutions more quickly. 
We refer to this as \emph{restricted SCOTS} (RSCOTS), where the restriction is in the number of switching actions allowed, or in the set of switchable arcs.

The available topologies are restricted in two ways.
First, the set of switchable arcs is reduced to $\wsarcs\subseteq\Arcs$.
 All other arcs are explicitly fixed to their default value:
\begin{equation}
	\z{u}{v} = \zDefault{u}{v} \quad \forall (u,v) \in \Arcs\setminus\wsarcs.
\end{equation} 
One option for selecting $\wsarcs$ is to rank arcs by a combination of the relevance scores from \Cref{subsec:partial-decomposition}, with weights $\wsWcp, \wsWsns$ and $\wsWvio \in \R$.
The number of switchable arcs is again limited by an absolute bound $\wsmaxn \in \N$, a relative bound $\wsmaxrel \in [0,1]$, and a score threshold $\wsminscore \in [0,1]$.
Our computational study instead uses the LODF-corrected susceptance sensitivity score defined below.
We first include arcs that are switched off in the default topology, then select the remaining arcs by decreasing score, up to a total of $\wsmaxn$ arcs.

Additionally, we reduce the global switching-action limit $\TAL$ (cf.\ Constraint \eqref{eq:switch-lim-all}).
The remaining constraints and objective are unchanged.
Note that reducing $\TAL$ to $0$ or setting $\wsarcs = \emptyset$ yields a simple SCOPF problem, without switching.
The sequential switching heuristic starts from SCOPF and solves restricted problems with switching limits $1,\ldots,\TAL$, passing on the best feasible solution between stages.
After accepted stages, it updates the arc scores and retains previously switched arcs as candidates.
It can be used as a standalone heuristic or to provide a warmstart for the exact solve.
Alternatively, the exact SCOTS solve runs with a parallel process executing the sequential heuristic.
The two processes exchange feasible incumbent solutions, so the main solve can benefit from the restricted problems without waiting for them to finish.

\paragraph{LODF-corrected susceptance sensitivity}
\label{par:susceptance-sensitivity}
For each arc $a \in \Arcs$ switched on in the default topology, let $\FlowCOptSym_a \in \R$ denote its SCOPF flow and $\DualPECOptSym_a \in \R$ the combined dual value of Constraints~\eqref{eq:sp-ohm-lb}--\eqref{eq:sp-ohm-ub} for contingency $\C \in \aC$.
We recover dual values from the MP solved as an LP with discrete decisions fixed, weighting contingency SP dual values by the MP dual multipliers of their Benders cuts.
The score is
\begin{equation}
  \wsSadm{a}
  = \frac{1}{1-h_a} \sum_{\C\in\aC}
        \FlowCOptSym_a\DualPECOptSym_a,
  \label{eq:susceptance-sensitivity}
\end{equation}
where $h_a$ is the change in flow on $a=(u,v)$ caused by a unit power transfer from $u$ to $v$ in the default topology.
Arcs with $1-h_a = 0$ receive score zero.
If $\FlowCOptSym_a\DualPECOptSym_a$ is positive, a small decrease in the arc's susceptance is predicted to reduce the objective value.
The factor $1/(1-h_a)$ accounts for the redistribution of power flows when arc $a$ is removed, so the score estimates the objective reduction from switching it off while keeping generation fixed.

\section{Evaluation}\label{sec:evaluation}
In this section we assess our branch-and-Benders-cut algorithm on a benchmark of 130 instances. 
Extensive experiments for parameter selection for all the algorithmic components were performed on a separate set of 31 instances, with most parts of the methodology being tuned sequentially.
Cut filtering and subproblem selection as well as partial decomposition and the sequential heuristic were tuned jointly.

\paragraph*{Test sets}
All benchmark instances were constructed from instances in the literature using the planar instance generation methodology proposed by~\textcite{donkiewicz2026instances}.
Instances of the PGLib-OPF~\cite{babaeinejad2021pglib} are selected randomly, reduced to smaller instances using a randomized breadth-first search, and then projected to 2D using a Kamada--Kawai layout algorithm~\cite{kamada1989algorithm}, before being connected with coupling edges.
Since these synthetic instances do not contain data on substations, generator groups, transformers, symmetry requirements, or emergency flexibility, all results are reported without the corresponding constraints.
For the evaluation on real-world instances including substation reconfiguration, we refer to the operational evidence in \Cref{sec:practice}.

\paragraph*{Implementation and setup}
We conducted all experiments on compute nodes running Debian~11 with two Intel Xeon L5630 processors at \qty{2.13}{\giga\hertz} and 16 hardware threads per node, with \qty{128}{\giga\byte} memory, four Gurobi threads for the master MILP, and one thread for the Benders and subproblem computations.
The implementation uses Julia~1.12.6~\cite{bezanson2017julia}, JuMP~1.30.1~\cite{lubin2023jump}, and PowerModels~0.21.6~\cite{coffrin2018powermodels}, with Gurobi~13.0.2~\cite{gurobi2026} across the experiments.

The sequential switching heuristic first solves SCOPF and then a sequence of RSCOTS problems, starting with 16 switchable arcs.
When used as a warmstart or parallel heuristic for the exact solve, it uses the same selection of switchable arcs:
we select initially open arcs first and the remaining arcs according to the LODF-corrected susceptance sensitivity score~\eqref{eq:susceptance-sensitivity}.
Unless stated otherwise, we include warm-starting in the reported solve times; for the exact solve with a parallel heuristic, we report the main solve time.
With a parallel heuristic, the four threads are split evenly between the restricted and main solve.

\subsection{Methodological Evaluation}
For the methodological evaluation, we fix the switching limit to $\TAL=1$, and evaluate different algorithmic components on a set of 130 instances.
We impose a time limit of \qty[drop-zero-decimal=true]{3}{\hour}, including warmstarting where applicable, and limit both SCOPF and each restricted warmstart stage to \qty{30}{\minute}.

\begin{table}[H]
  \centering\scriptsize
  \setlength{\tabcolsep}{2.3pt}
  \caption{Comparison of selected algorithm configurations on 130 SCOTS instances.}
  \label{tab:final-parameter-study}
  \begin{tabular}{@{}rcl>{\raggedright\arraybackslash}p{0.40\linewidth}r@{\hspace{12pt}}rr@{}}
\toprule
 & Pred. & Method & Configuration & Solved & Mean (s) & Geo. Mean (s) \\
\midrule
1 & - & Compact & All contingencies in master & 110 & 2471.20 & 534.37 \\
2 & 1 & Benders & LP subproblems & 129 & 921.91 & 403.46 \\
3 & 2 & Benders SLE-U & SLE-U subproblem solver & 130 & 140.73 & 92.09 \\
4a & 3 & SP Selection & Random order; $\spmaxn=2$, $\spstab=100$ & 130 & 72.04 & 59.71 \\
4b & 3 & Cut Filtering & Violation ranking; $\cfmaxn=10$ & 130 & 92.66 & 72.59 \\
5 & 4 & SP+CF combined & Replace $\spmaxn=2$ by $\spmaxrel=0.05$; $\cfmaxn=5$; no stabilization rounds & 130 & 64.96 & 57.13 \\
6 & 5 & Partial Decomp. & Base-case flow in master; $\pdmaxn=0$ & 130 & 57.16 & 52.34 \\
7 & 6 & Branch Fixing & Fix branches required for connectivity to on & 130 & 52.07 & 48.11 \\
\textcolor{gray}{8a} & \textcolor{gray}{7} & \textcolor{gray}{Heuristic warmstart} & \textcolor{gray}{SCOPF, then RSCOTS with 16 switchable arcs} & \textcolor{gray}{130} & \textcolor{gray}{56.49} & \textcolor{gray}{51.61} \\
\textcolor{gray}{8b} & \textcolor{gray}{7} & \textcolor{gray}{Parallel heuristic} & \textcolor{gray}{SCOPF--RSCOTS concurrent with SCOTS; two threads per process} & \textcolor{gray}{130} & \textcolor{gray}{54.24} & \textcolor{gray}{50.21} \\
\textcolor{gray}{9} & \textcolor{gray}{7} & \textcolor{gray}{Multithreading} & \textcolor{gray}{4 Benders subproblem threads} & \textcolor{gray}{130} & \textcolor{gray}{39.05} & \textcolor{gray}{36.86} \\
\bottomrule
\end{tabular}
% n = 130; excluded pglib_scopf_case162; shift = 10 s; non-OPTIMAL/invalid = 10800 s.

\end{table}
\begin{figure}[H]
  \centering
  \begin{subfigure}{\linewidth}\centering
    \pgfplotsset{every axis post/.append style={height=0.295\linewidth}}
    \input{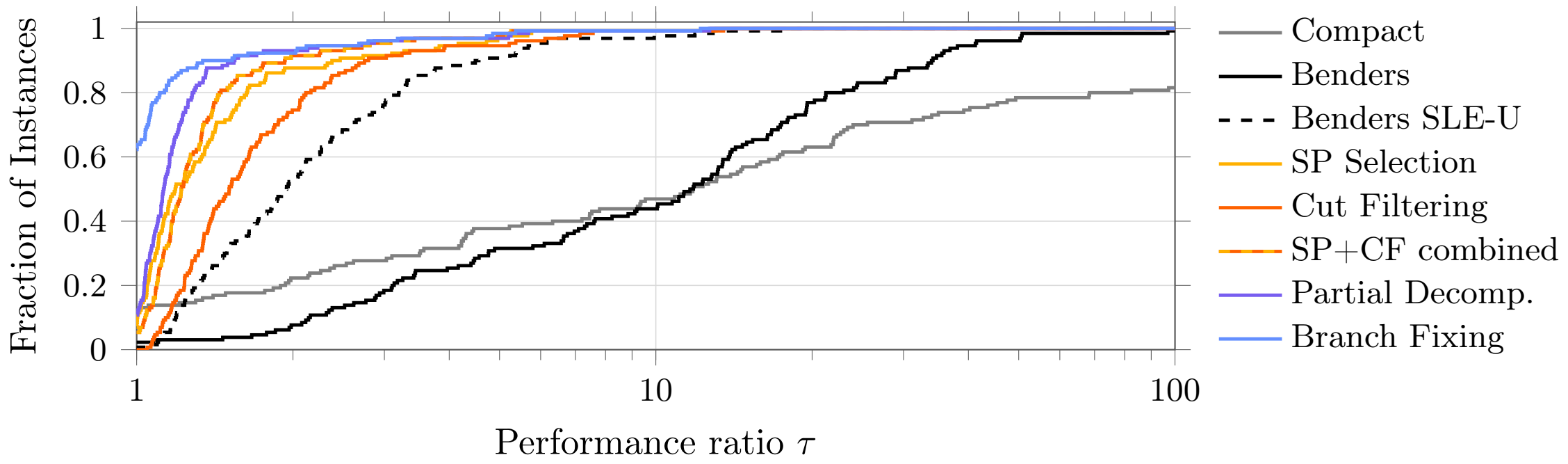}
  \end{subfigure}
  \caption{Performance profile of the configurations in \Cref{tab:final-parameter-study}, excluding grey rows.}
  \label{fig:final-pp}
\end{figure}

We summarize the final settings and the performance in \Cref{tab:final-parameter-study} and visualize the total solve times in a performance profile in \Cref{fig:final-pp}.
``Geo. Mean'' is the shifted geometric mean with a \qty{10}{\second} shift.
Both means assign \qty{10800}{\second} to unsuccessful or invalid runs.

The final configuration (``Branch Fixing'') solves all 130 instances, compared with 129 for Benders decomposition with LP subproblems (and 110 for the compact formulation), and reduces the shifted geometric mean solve time from \qty{403.46}{\second} to \qty{48.11}{\second}, a reduction of \qty{88.1}{\percent}.
It includes all discussed components, which address different computational bottlenecks.
Replacing LP subproblems by SLE-U increases the share of solve time attributed to master problems from \qty{7.6}{\percent} to \qty{51.9}{\percent} on the 129 instances solved by both variants.
With master problem solving now being the main bottleneck, it is crucial to control its growth.
To that end, subproblem selection limits the number of subproblems evaluated (and thus also cuts added), while cut filtering directly controls the number of cuts.
With jointly tuned limits of subproblem selection and cut filtering, their combination reduces the shifted geometric mean to \qty{57.13}{\second}, compared with \qty{59.71}{\second} for selection alone and \qty{72.59}{\second} for filtering alone.
Partial decomposition with only the base-case power flow in the master yields a further \qty{8.4}{\percent} reduction to \qty{52.34}{\second}.
One can also add contingency power flows to the master instead of the base case, but we were not able to observe a statistically significant advantage in doing that.
Fixing branches required for connectivity to on reduces the shifted geometric mean solve time by another \qty{8.1}{\percent}, to \qty{48.11}{\second}.

Using four Benders subproblem threads reduces the shifted geometric mean solve time by a further \qty{23.4}{\percent}, to \qty{36.86}{\second}, but it should be noted that this is an uneven comparison.
For $Z=1$, heuristic warmstarting and a parallel heuristic provide no substantial runtime benefit.
We omit the heuristic warmstart, parallel heuristic, and multithreading configurations from the performance profile in \Cref{fig:final-pp}.

\subsection{Practical Evaluation}
In this subsection, we evaluate our method in a more practical setting, with higher switching limits $\TAL$ and a strong focus on primal convergence.
We compare the exact approach with only running the sequential switching heuristic or a single-stage variant of it (solving the last stage of the sequential heuristic, but for the full time budget).
Additionally, we compare with two combinations of the exact approach and the heuristic: either warmstarting the exact solve with the heuristic solution or running the heuristic in parallel with the exact solve.
We run each configuration separately for switching limits $Z=1,\ldots,5$, with a six-hour budget per run.
For the upper bounds at $Z=3$ in \Cref{fig:cumulative-savings}, we additionally use runs with a \qty{96}{\hour} time limit.
SCOPF and restricted stages are capped at \qty{30}{\minute}, except the parallel heuristic's final stage, which may continue until the main solve or budget ends.

We require an absolute objective reduction of at least \num[round-precision=4]{0.0001} to count an improvement to reduce numerical noise.

\newlength{\practicalFigureWidth}
\setlength{\practicalFigureWidth}{\linewidth}

\paragraph*{Sequential switching heuristic}
The single-stage and sequential switching heuristics use the same initial 16 candidates selected from the SCOPF solution and total allowance of $\num{1800} \cdot \TAL$ seconds (\qty{1800}{\second} per sequential stage), plus \qty{1800}{\second} for SCOPF.
Both use zero relative MIP gap without objective-based stopping or a subsequent exact solve.
The single-stage heuristic keeps its candidate set fixed.
The sequential heuristic updates the ranking and retains previously switched branches, so its candidate set may contain more than 16 branches for some stages after the first.

\begin{figure}[H]
  \centering
  \pgfplotsset{every axis post/.append style={axis on top}}
  \input{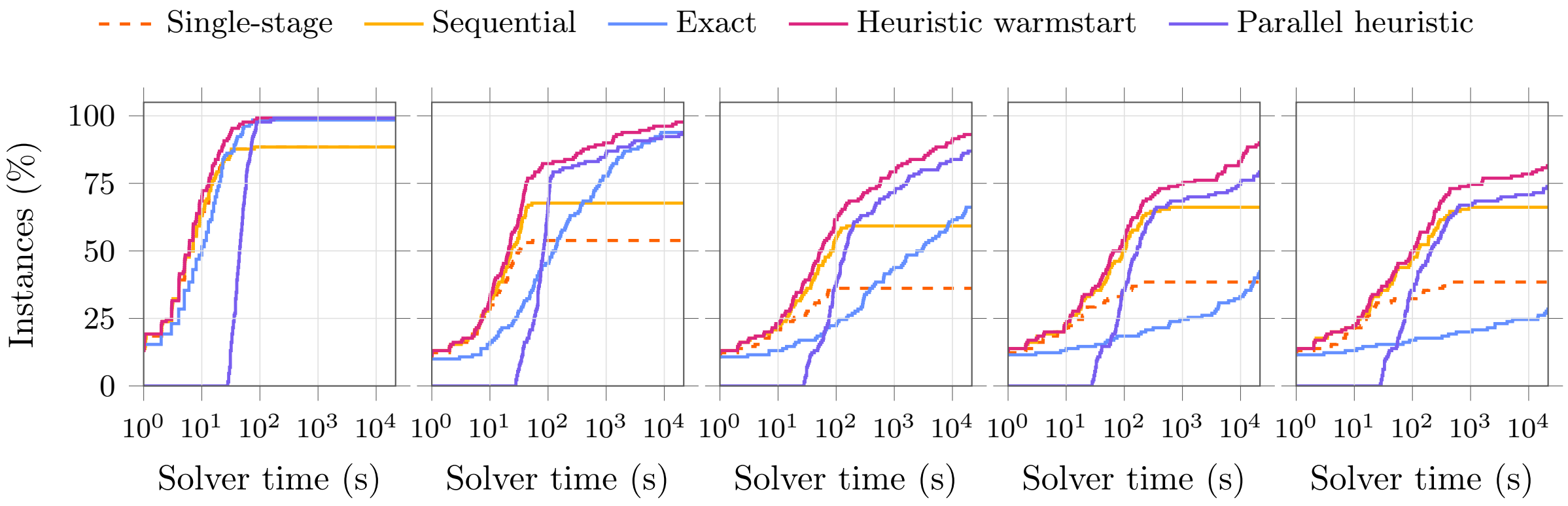}
  \caption{Instances within \num[round-precision=4]{0.0001} of the best objective among the five configurations, excluding SCOPF time.}
  \label{fig:heuristic-selection}
\end{figure}

\Cref{fig:heuristic-selection} shows the fraction of instances where each configuration finds a solution that is within \qty[round-precision=4]{0.0001} of the best known objective among all five configurations.
One can see that between the sequential heuristic and the single-stage heuristic, the former finds the best known solution significantly more often for $Z\geq2$.
For $Z\leq3$, the exact solve alone is able to find the best known solution on a larger fraction of instances than the sequential heuristic, but with significantly longer running time.
For $Z\geq4$, it finds the best known solution less frequently than the sequential heuristic.
For all $Z$, both the heuristic warmstart and the parallel heuristic find the best known solution as often or more often than the other three configurations, with warmstart reaching it more often for $Z\geq2$.
This may be due to the parallel heuristic's threads being split between the heuristic and exact solve, which may reduce the effectiveness of both.
We conclude that while the heuristics perform reasonably well, the exact solve is still sensible for finding even better solutions, especially for larger switching limits.

We focus only on heuristic warmstarting in the following.

\begin{figure}[H]
  \centering
  \begin{subfigure}[t]{0.385\linewidth}
    \centering
    % Absolute reduction threshold: 0.0001; six-hour cutoff before SCOPF removal.
% Z=3 bounds include additional runs with up to 96 h wall time (95 h solver time).
% Solid curves retain the original six-hour solutions.
\def\cumulativeSavingsN{130}
\ifdefined\practicalFigureWidth\else
\newlength{\practicalFigureWidth}\setlength{\practicalFigureWidth}{\linewidth}\fi
\begin{tikzpicture}[baseline=(current axis.north)]
\definecolor{practicalSequential}{HTML}{DC267F}
\begin{axis}[scale only axis,width=0.255\practicalFigureWidth,height=0.18\practicalFigureWidth,
xmin=0,xmax=5,ymin=0,ymax=30,xtick={0,1,2,3,4,5},ytick distance=5,
xlabel={Switching limit $Z$},ylabel={Objective reduction\\relative to SCOPF (\%)},
ylabel style={align=center},
grid=major,major grid style={draw=black!12,line width=0.3pt},xmajorgrids=false,
axis line style={draw=black!65,line width=0.4pt},tick align=outside,
tick style={draw=black!65},label style={font=\small},tick label style={font=\footnotesize},
legend columns=2,legend cell align={left},
legend style={at={(0.5,1.08)},anchor=south,draw=none,fill=none,font=\footnotesize,
/tikz/every even column/.append style={column sep=0.9em}},]
\addplot[draw=none,fill=practicalSequential!14,forget plot] coordinates {(0,0) (1,11.8486862) (2,15.07034278) (3,15.81343755) (4,15.96186107) (5,15.9796788) (5,27.79647489) (4,24.4201359) (3,17.01848474) (2,15.09276631) (1,11.84923827) (0,0)} \closedcycle;
\addplot[draw=none,fill=practicalSequential!7,forget plot] coordinates {(0,0) (1,2.216687) (2,3.655361948) (3,3.885294365) (4,3.664821749) (5,3.670263042) (5,11.4472084) (4,8.680612972) (3,4.430548811) (2,3.659922821) (1,2.216687) (0,0)} \closedcycle;
\addplot[color=practicalSequential,dashed,no marks,line width=0.7pt,forget plot] coordinates {(0,0) (1,11.84923827) (2,15.09276631) (3,17.01848474) (4,24.4201359) (5,27.79647489)};
\addplot[color=practicalSequential,solid,mark=*,mark options={fill=practicalSequential},line width=0.9pt,mark size=1.8pt] coordinates {(0,0) (1,11.8486862) (2,15.07034278) (3,15.81343755) (4,15.96186107) (5,15.9796788)};
\addlegendentry{Mean}
\addplot[color=practicalSequential,dashed,no marks,line width=0.7pt,forget plot] coordinates {(0,0) (1,2.216687) (2,3.659922821) (3,4.430548811) (4,8.680612972) (5,11.4472084)};
\addplot[color=practicalSequential,solid,mark=square*,mark options={fill=white},line width=0.9pt,mark size=1.8pt] coordinates {(0,0) (1,2.216687) (2,3.655361948) (3,3.885294365) (4,3.664821749) (5,3.670263042)};
\addlegendentry{Median}
\end{axis}
\end{tikzpicture}
    \edef\practicalCaption{Objective reductions ($n=\cumulativeSavingsN$). Dashed curves give upper bounds.}
    \expandafter\caption\expandafter{\practicalCaption}
    \label{fig:cumulative-savings}
  \end{subfigure}\hfill
  \begin{subfigure}[t]{0.605\linewidth}
    \centering
    \input{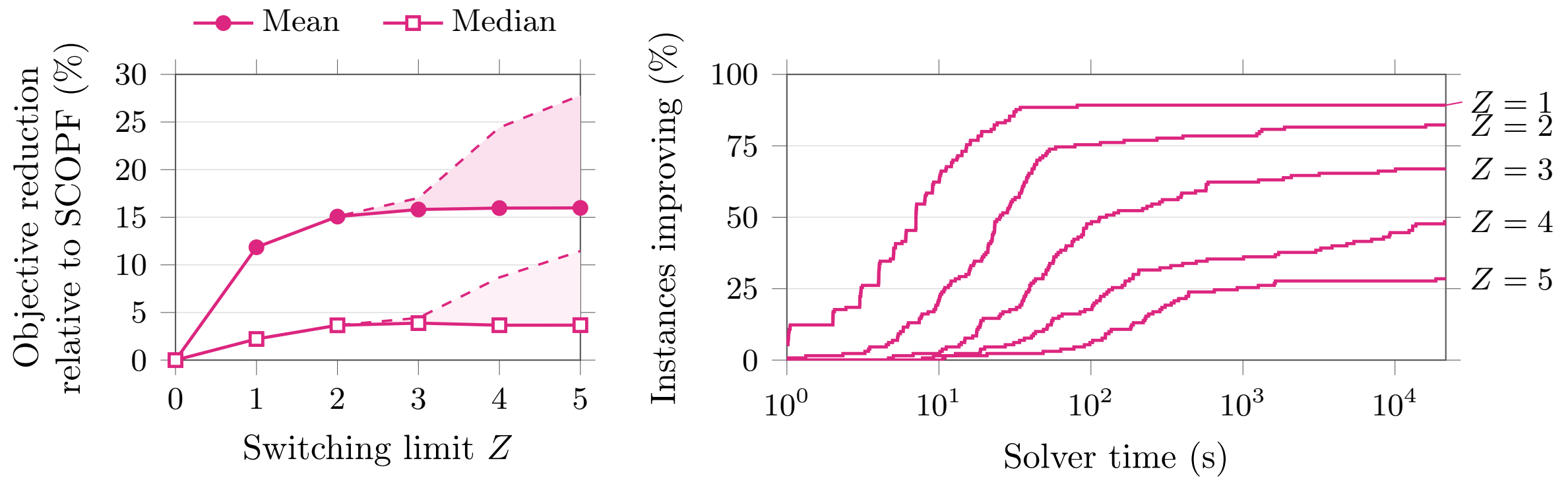}
    \caption{Time to improve on the best solution at $Z-1$ ($n=130$), excluding SCOPF time.}
    \label{fig:stepwise-improvements}
  \end{subfigure}
  \caption{Solution quality and improvement times for the exact solve with heuristic warmstart.}
  \label{fig:practical-quality-time}
\end{figure}

\paragraph*{Objective reductions}
Most of the observed mean reduction is achieved by $Z=2$ (\Cref{fig:cumulative-savings}), while the much smaller median indicates substantial variation across instances.
The widening gaps to the upper bounds leave the potential benefit of further switching actions unresolved.
Runs at higher $Z$ may fail to find the best solution for smaller $Z$ within the time limit because the problem tends to get harder with more switching actions, and even with sequential heuristic as warmstart, not all optimal solutions of the previous switching limit are found.

\paragraph*{Improvement times}
\Cref{fig:stepwise-improvements} shows the percentage of instances for which an improvement over the best known solution with one fewer switching action is found over time.
Most notably, the percentage of instances that achieve improvements decreases with increasing $Z$, and the median time to find such an improvement increases.
For all $Z$, most successfully improving runs do so within \qty{10}{\minute}: 31 of the 37 improvements at $Z=5$ occur within this interval, and six occur later.
Additional switching actions therefore remain useful on some instances, while finding an improvement becomes more difficult within the available time.

\section{Practical Application}\label{sec:practice}
We finally discuss additional insights stemming from our application of the introduced algorithm on real-world data in collaboration with the German TSOs.
We comment on the use of the algorithm in an applied setting and its practical benefits and limitations.

\subsection{Practical Performance and Applicability}
\label{subsubsec:feasibility-practice}
A previous version of the algorithm (similar to configuration 5 in \autoref{tab:final-parameter-study}) was evaluated on instances representing certain historical situations in the German grid, including some nodes and branches in neighboring countries, with flows over branches going further out modeled as static.
Different instances use an almost identical network, consisting of roughly 3000 buses, 5000 branches, 2300 generators, and 2800 loads.
In contrast to our evaluation, where all individual branches that do not trivially split the network are considered switchable and as possible contingencies, there are predefined subsets of branches for both.
This is because some branches do not represent switchable components, and components in the data do not necessarily map exactly to real-world components. 
In the final model, over 300 substations of sizes ranging from 2 to 5 busbars were considered, along with over 1300 switchable branches and roughly 1000 contingencies.

\paragraph*{Results on real-world data}
A first evaluation was conducted on 71 historical operating days, with 24 hourly instances per day.
The instances allowed up to five line openings and no substation switching.
For all of them, a mathematically feasible topology was found within a five-hour time limit~\cite{innoptem-abschlussbericht}.
Relative to the solution without topology optimization, our algorithm achieved a reduction of approximately \qty{7}{\percent} in simulated redispatch volume, with a median total CPU time of approximately \qty{34}{\minute} per instance.
Since data for the actual topological actions executed by the TSOs would not represent a fair comparison (during the day, more accurate forecasts may become available and enable better choices of switching actions, and our model omits some operational constraints), we cannot compare our actions to the actual actions that TSOs are already taking.
These results therefore do not establish additional redispatch cost savings over current practice.
It should be noted that significant computational improvements were achieved between the algorithm used for the real-world data and its state as reported for the synthetic instances.

\paragraph*{Real-world feasibility}
Next to the performance of the derived instances, it is important to assess practical feasibility, as there may be additional operational constraints that are not captured by our model.
In particular, apart from potentially unknown or omitted constraints, the DC approximation in our model does not necessarily yield AC-feasible solutions.
Feasibility of selected solutions was assessed through a series of workshops conducted jointly with all four German TSOs.
We presented solutions of seven different instances generated by our algorithm to the attendant operating personnel.
The solutions included both line switching and substation reconfiguration, and upon initial assessment, all of them were considered operationally plausible, though further assessment would be necessary.
It has to be noted though that operators emphasized that some of the proposed actions may not always be worth the effort, \eg, to coordinate with power plant personnel. 
This also underscores the importance of maintaining a human-in-the-loop after the optimization phase.
All solutions passed the subsequent AC feasibility verification by the TSOs.

\subsection{Operational case study}
\label{subsubsec:case-study}

We illustrate the integrated approach on a forecasted operating hour from the project final report~\cite{innoptem-abschlussbericht}.
The network contains two congestion regions. In the larger region, the two most heavily loaded circuits reach \qty{132}{\percent} and \qty{130}{\percent} of their limits under a single outage. Without topology changes, resolving the congestion requires more than \qty{1}{\giga\watt\hour} of simulated redispatch.
The integrated optimization opens two busbar couplers, separating one substation into three busbars.
In the subsequent AC power flow assessment, this topology change reduces the two loadings to \qty{107}{\percent} and \qty{110}{\percent} before redispatch. 
According to simulations, the combination of topology changes and redispatch reduces the calculatory redispatch cost by \qty{20.4}{\percent} relative to redispatch without topology changes.

The voltage magnitude difference across the opened couplers reaches \qty{14.2}{\kilo\volt} under an outage and \qty{7.3}{\kilo\volt} in the base case, below the study limit of \qty{15}{\kilo\volt}.
The report nevertheless requires an operator assessment before implementation, taking account of the network situation and fulfillment of other constraints that are not yet modeled. 
Meeting the study limit alone therefore does not establish that the proposed action can be implemented.

\subsection{Key Benefits and Limitations of the Integrated Optimization}\label{subsubsec:benefits-limitations}
Due to the difficulty and size of the problem, certain simplifications can be taken in order to make it tractable (notably, we consider AC only after fixing the topology and optimize hours sequentially, which can yield suboptimal results).
We discuss several such simplifications and provide computational evidence that they can lead to suboptimal or even infeasible solutions.

\paragraph{Atomic switching action ranking}
A heuristic approach can be to consider switching actions atomically, \ie, multiple switching actions might be considered, but they are selected sequentially based on their individual impact on some network state.
This can miss sets of switching actions that are only effective when taken together.
An example of this is the decoupling of two strongly coupled sectors of the grid: if they are connected by two parallel lines (which are both overloaded), switching off one of them even increases the overload on the other line further.
Deactivating both lines, however, can relieve the congestion by rerouting the power flow through other paths.
In 18 of 107 optimal two-switch solutions on the synthetic test set (see \Cref{sec:evaluation}), neither switch occurs in the corresponding $\TAL=1$ solution.
For 28 of 107 optimal two-switch solutions on the synthetic test set, at least one switch individually increases cost even after redispatch is reoptimized.
We also observed this for real-world instances.

\paragraph{Pre-filtered contingencies}
Similarly to switching actions, contingencies can be filtered based on whether they are (close to) violated for some initial solution, with the assumption that any actions that are taken will not realistically lead to a violation. 
This is, however, not always the case, so pre-filtering without checking all contingencies for the final solution can admit infeasible solutions.
In one instance, a contingency omitted by another approach became critical after switching.
The resulting N\nobreakdash-1 violation was detected only in the final solution check.

\paragraph{Restricted set of topologies}
In practice, the high number of possible topologies poses a computational challenge, so some tools may restrict themselves to a comparatively small set of topological actions.
Extending the allowed region, for example, from mandating a maximum of two connected groups of busbars (i.e., two-node operation of a substation) to arbitrary busbar splits (in Germany, there are substations with up to five busbars) is a complex undertaking if there is no computational decision support.
The combinatorial explosion of possible topologies makes it hard to check for the extensive operational constraints and thus gain trust in the new topologies in the short timeframe operators have to approve the proposed actions.
For example, in one proposed solution, we found a multi-busbar split into three nodes (\emph{Dreiknotenbetrieb}) to be the optimal solution, which is a less common solution compared to the conventional busbar-split into two nodes.
Efficient computational support can facilitate adopting such unconventional but cost-effective solutions.

\paragraph{Sequential optimization of topology and redispatch}
A common pattern, also used in practice, is to sequentially optimize topology and redispatch, \ie, redispatch is optimized with a fixed topology, then the topology is optimized with that redispatch fixed and so on.
Like with sequential switching actions, this can miss redispatch and topology solutions that only work together.
For example, on real-world data, we were able to reduce simulated redispatch costs by shifting congestion to regions where it could be relieved more cheaply~\cite[Section~6.2.2]{innoptem-abschlussbericht}.
This topology change did not alleviate the congestion (it was provably not possible to do so at all for the instance), but reduced cost.
Additionally, redispatch may be used not only to reduce congestion, but also to adhere to constraints such as angle limits after switching, which may allow for more effective topology actions.
Sequential optimization can exclude such solutions, because relocating an overload need not directly reduce costs.

\paragraph{Minimizing redispatch volume instead of cost}
In Germany, it is mandated by law that calculatory cost is minimized in the preventive redispatch process.
A common misconception that we encounter is that redispatch volume can be equated with redispatch costs.
However, not only do differing redispatch costs already indicate that this may not always be the case, we were also able to empirically show that optimal solutions can even increase redispatch volume, while reducing calculatory costs.

\section{Concluding Remarks}\label{sec:summary}
We present a mixed-integer linear programming approach to solve a real-world SCOTS problem, where congestion in high-voltage power grids is managed using network and generator (``market'') measures across a high number of contingency scenarios.
A shortfall of current practice is the separated optimization of network and market measures.
We give a complete, real-world applicable MILP model, which is then solved with a branch-and-Benders cut algorithm enhanced by several relevant additions.
Its performance is tested on synthetic and real-world instances, and we demonstrate that it is able to find provably optimal solutions in reasonably small time.
Finally, we discuss potential issues with heuristic approaches.

Upon the acceptance of this paper, the software that implements the described algorithm, including several data processing steps not described in the paper, will be made publicly available as open-source.
To our knowledge, there exists no such software package that implements a complete, realistic model and provides the methodology to solve it to provable optimality.
Due to the nature of MILP, the model can be easily modified by adding or changing the constraints, variables, and objective function, without a need to alter the algorithm itself, as those changes will be enforced automatically.
Since we propose no sequential solving (the sequential heuristic reuses the same model) or manual checking of constraints, this does not require understanding the algorithm, or even most code. 
This enables TSOs (also outside Germany), regulators, and individuals to optimize the problem and also conduct sensitivity studies around it.
The presented method can also be transferred to curative measures (i.e., reactive switching and redispatch), and even serve use in transmission expansion planning without the need to completely starting from scratch.

Directions for future work are thus numerous and not constrained to further improving the algorithm, but also include broader applications of optimization in the energy sector.

\par\smallskip\noindent\textbf{Acknowledgements.}\enspace
We thank Felix Preuschoff for his electrical engineering expertise and support in developing the full operational model.

\par\smallskip\noindent\textbf{Author contributions: CRediT.}\enspace
\textbf{Tim Donkiewicz:} Conceptualization, Data curation, Formal analysis, Investigation, Methodology, Software, Validation, Visualization, Writing -- original draft, Writing -- review \& editing.
\par\noindent\textbf{Oliver Gaul:} Conceptualization, Formal analysis, Investigation, Methodology, Software, Validation, Visualization, Writing -- original draft, Writing -- review \& editing.
\par\noindent\textbf{Marco L\"ubbecke:} Conceptualization, Funding acquisition, Methodology, Supervision, Writing -- review \& editing.

\par\smallskip\noindent\textbf{Data and code availability.}\enspace
The code, public test instances, and reproduction scripts are available in \href{https://github.com/tidonk/TopologyOptimizer.jl/tree/integrated-2026-11-10}{TopologyOptimizer.jl (tag \texttt{integrated-2026-11-10})}.
The operational network data are confidential and cannot be shared.

\par\smallskip\noindent\textbf{Declaration of competing interests.}\enspace
The authors declare that they have no competing interests.

\par\smallskip\noindent\textbf{Declaration of generative AI and AI-assisted technologies in the manuscript preparation process.}\enspace
The authors used large language model tools, including GitHub Copilot and OpenAI Codex, for assistance with code, literature searches, manuscript editing, and notation checks.
They reviewed and edited the content and take full responsibility for the published article.

%% bibliography %%
\printbibliography

\end{document}